\documentclass[12pt]{article} 

\usepackage[dvips]{graphics,lscape}
\DeclareGraphicsExtensions{.eps.gz,.eps,.epsi.gz,.epsi,.ps,.ps.gz}
\usepackage{epsfig}
\usepackage{color}
\graphicspath{{fig/}} \textwidth=6in \textheight=9in

\usepackage{latexsym}
\usepackage{graphicx}
\usepackage{amsfonts}
\usepackage{amssymb}
\usepackage{mathrsfs}

\usepackage{amsmath}
\usepackage{latexsym}
\usepackage{mathrsfs}
\usepackage{amsfonts}
\usepackage{amssymb}

\def\argmin{\mathop{\rm arg\, min}}

\newcommand{\bel}{\begin{eqnarray}\label}
\newcommand{\eel}{\end{eqnarray}}
\newcommand{\bes}{\begin{eqnarray*}}
\newcommand{\ees}{\end{eqnarray*}}

\def\benu{\begin{enumerate}}
\def\eenu{\end{enumerate}}

\def\bhbar{{\overline \bh}}
\def\teps{{\widetilde\eps}}

\def\argmin{\mathop{\rm arg\, min}}
\def\real{{\mathbb{R}}}
\def\R{{\real}}

\def\complex{\mathop{{\rm I}\kern-.58em\hbox{\rm C}}\nolimits}
\def\pa{\partial}

\def\diag{\hbox{\rm diag}}

\def\sgn{\hbox{\rm sgn}}
\def\trace{\hbox{\rm trace}}

\def\Var{\hbox{\rm Var}}

\def\supp{\hbox{\rm supp}}

\def\mathbold{\boldsymbol} 

\def\bA{\mathbold{A}}\def\Atil{{\widetilde A}}

\def\bb{\mathbold{b}}
\def\hbb{{\widehat{\bb}}}
\def\bB{\mathbold{B}}

\def\bD{\mathbold{D}}\def\Dhat{\widehat{D}}
\def\hbD{{\widehat{\bD}}}

\def\bfe{\mathbold{e}}

\def\bg{\mathbold{g}}

\def\bh{\mathbold{h}}
\def\tbh{{\widetilde{\bh}}}

\def\bI{\mathbold{I}}

\def\Ltil{{\widetilde L}}

\def\bM{\mathbold{M}}

\def\bP{\mathbold{P}}

\def\bR{\mathbold{R}}

\def\Shat{\widehat{S}}

\def\bu{\mathbold{u}}

\def\bU{\mathbold{U}}
\def\scrU{{\mathscr U}}

\def\bv{\mathbold{v}}

\def\bx{\mathbold{x}}

\def\bX{\mathbold{X}}

\def\by{\mathbold{y}}

\def\bz{\mathbold{z}}

\def\balpha{\mathbold{\alpha}}

\def\hbalpha{{\widehat{\balpha}}}

\def\bbeta{\mathbold{\beta}}\def\hbeta{\widehat{\beta}}

\def\hbbeta{{\widehat{\bbeta}}}\def\tbbeta{{\widetilde{\bbeta}}}
\def\bbetabar{{\overline\bbeta}}

\def\ep{\varepsilon}\def\eps{\epsilon}
\def\bep{\mathbold{\ep}}

\def\bTheta{\mathbold{\Theta}}\def\hbTheta{{\widehat{\bTheta}}}
\def\tbTheta{{\widetilde{\bTheta}}}

\def\Thetatil{{\widetilde \Theta}}

\def\lam{\lambda}
\def\hlam{\widehat{\lam}}
\def\tlam{\widetilde{\lam}}

\def\bsigma{\mathbold{\sigma}}\def\hsigma{\widehat{\sigma}}
\def\sigmabar{{\overline{\sigma}}}
\def\hbsigma{{\widehat{\bsigma}}}
\def\bsigmabar{{\overline\bsigma}}

\def\bSigma{\mathbold{\Sigma}}

\def\Sigmabar{{\overline \Sigma}}
\def\bSigmabar{{\overline\bSigma}}

\def\htau{\widehat{\tau}}

\def\hphi{\widehat{\phi}}

\def\bOmega{\mathbold{\Omega}}\def\hbOmega{{\widehat{\bOmega}}}
\def\tbOmega{{\widetilde{\bOmega}}}

\newtheorem{theorem}{Theorem}
\newtheorem{lemma}{Lemma}
\newtheorem{proposition}{Proposition}
\newtheorem{corollary}{Corollary}
\newtheorem{example}{Example}
\newtheorem{remark}{Remark}

\def\argmin{\mathop{\rm arg\, min}}

\def\teps{{\widetilde\eps}}

\def\bSigmabar{\mathbold{\Sigmabar}}
\def\bhbar{{\overline \bh}}

\begin{document}
\title{Sparse Matrix Inversion with Scaled Lasso}
\author{
Tingni Sun\\
Statistics Department, The Wharton School, University of Pennsylvania\\ Philadelphia, Pennsylvania, 19104\\
{tingni@wharton.upenn.edu} \\
\\
Cun-Hui Zhang\\
Department of Statistics and Biostatistics, Rutgers University\\
Piscataway, New Jersey 08854\\
{czhang@stat.rutgers.edu}
}
\date{}
\maketitle 

\bigskip\begin{abstract}
We propose a new method of learning a sparse nonnegative-definite target matrix. 
Our primary example of the target matrix is the inverse of a population covariance or correlation matrix. 
The algorithm first estimates each column of the target matrix by the scaled Lasso 
and then adjusts the matrix estimator to be symmetric. 
The penalty level of the scaled Lasso for each column is completely determined 
by data via convex minimization, without using cross-validation.

We prove that this scaled Lasso method guarantees the fastest proven rate of convergence 
in the spectrum norm under conditions of weaker form than those in the existing analyses of 
other $\ell_1$ regularized algorithms, and has faster guaranteed rate of convergence when 
the ratio of the $\ell_1$ and spectrum norms of the target inverse matrix diverges to infinity.
A simulation study demonstrates the computational feasibility and superb 
performance of the proposed method.

Our analysis also provides new performance bounds for the Lasso and scaled Lasso 
to guarantee higher concentration of the error at a smaller threshold level than previous analyses, 
and to allow the use of the union bound in column-by-column applications of the scaled Lasso 
without an adjustment of the penalty level. 
In addition, the least squares estimation after the scaled Lasso selection is considered and proven 
to guarantee performance bounds similar to that of the scaled Lasso.
\end{abstract}

\section{{Introduction}}
We consider the estimation of the matrix inversion $\bTheta^*$ satisfying $\bSigmabar\bTheta^*\approx\bI$ 
for a given data matrix $\bSigmabar$. 
When $\bSigmabar$ is a sample covariance matrix, our problem is the estimation of the inverse of 
the corresponding population covariance matrix. 
The inverse covariance matrix is also called precision matrix or concentration matrix. 
With the dramatic advances in technology, the number of variables $p$, 
or the size of the matrix $\bTheta^*$, is often of greater order than 
the sample size $n$ in statistical and engineering applications. 
In such cases, the sample covariance matrix is always singular and 
a certain type of sparsity condition is typically imposed 
for proper estimation of the precision matrix and for theoretical investigation of the problem. 
In a simple version of our theory, this condition is expressed as the $\ell_0$ sparsity, 
or equivalently the maximum degree, of the target inverse matrix $\bTheta^*$. 
A weaker condition of capped $\ell_1$ sparsity is also studied to allow many small signals.

Several approaches have been proposed to the estimation of sparse inverse matrices in high-dimensional setting. 
The $\ell_1$ penalization is one of the most popular methods. Lasso-type methods, or convex minimization algorithms with the $\ell_1$ penalty on all entries of $\bTheta^*$, have been developed in \cite{BanerjeeEA08,FriedmanHT08}, and in \cite{YuanL07} with $\ell_1$ penalization on the off-diagonal matrix only. This is refereed to as the graphical Lasso (GLasso) due to the connection of the precision matrix to Gaussian Markov graphical models. 
In this GLasso framework, \cite{RavikumarWRY08} provides sufficient conditions for model 
selection consistency, while \cite{RothmanBLZ08} provides the convergence rate 
$\{((p+s)/n)\log p\}^{1/2}$ in the Frobenius norm and $\{(s/n)\log p\}^{1/2}$ in the spectrum norm, 
where $s$ is the number of nonzero off-diagonal entries in the precision matrix.  
Concave penalty has been studied to reduce the bias of the GLasso \cite{LamF09}.
Similar convergence rates have been studied under the Frobenius norm in a unified framework 
for penalized estimation in \cite{NegahbanRWY10}.
Since the spectrum norm can be controlled via the Frobenius norm, 
this provides a sufficient condition $(s/n)\log p \to 0$ for the convergence 
to the unknown precision matrix under the spectrum norm. 
However, in the case of $p\ge n$, this condition does not hold for 
banded precision matrices, 
where $s$ is of the order of the product of $p$ and the width of the band. 

A potentially faster rate $d\sqrt{(\log p)/n}$ 
can be achieved by $\ell_1$ regularized estimation of individual columns of the precision matrix, 
where $d$, the matrix degree, is the largest number of nonzero entries in a column. 
This was done in \cite{Yuan10} by applying the Dantzig selector to the regression of each variable 
against others, followed by a symmetrization step via linear programming. 
When the $\ell_1$ operator norm of the precision matrix is bounded, this method 
achieves the convergence rate $d\sqrt{(\log p)/n}$ in $\ell_q$ matrix operator norms. 
The CLIME estimator \cite{CaiLL11}, which uses the Dantzig selector directly to 
estimate each column of the precision matrix, also achieves the $d\sqrt{(\log p)/n}$ rate 
under the boundedness assumption of the $\ell_1$ operator norm. 
In \cite{YangK10}, the Lasso is applied to estimate the columns of the target matrix under the assumption of equal diagonal, and the estimation error is studied in the Frobenius norm for $p=n^{\nu}$.
This column-by-column idea reduces a graphical model to $p$ regression models. It was first introduced by   \cite{MeinshausenB06} for identifying nonzero variables in a graphical model, called neighborhood selection.
In addition, \cite{RochaZY08} proposed a pseudo-likelihood method by merging all $p$
linear regressions into a single least squares problem.

In this paper, we propose to apply the scaled Lasso \cite{SunZ11} column-by-column 
to estimate a precision matrix in the high dimensional setting.  
Based on the connection of precision matrix estimation to linear regression, 
we construct a column estimator with the scaled Lasso, 
a joint estimator for the regression coefficients and noise level.
Since we only need a sample covariance matrix as input, this estimator could be extended 
to generate an approximate inverse of a nonnegative-definite data matrix in a more general setting. 
This scaled Lasso algorithm provides a fully specified map from the space of nonnegative-definite matrices 
to the space of symmetric matrices. For each column, the penalty level of the scaled Lasso is determined by
data via convex minimization, without using cross-validation.

We study theoretical properties of the proposed estimator for a precision matrix under a normality assumption.
More precisely, we assume that the data matrix is the sample covariance matrix $\bSigmabar=\bX^T\bX/n$, 
where the rows of $\bX$ are iid $N(0,\bSigma^*)$ vectors. 
Let $\bR^* = (\diag\bSigma^*)^{-1/2}\bSigma^*(\diag\bSigma^*)^{-1/2}$ be the population 
correlation matrix. Our target is to estimate the inverse matrices 
$\bTheta^*=(\bSigma^*)^{-1}$ and $\bOmega^* = (\bR^*)^{-1}$. Define 
\bel{degree}
d = \max_{1\le j\le p}\#\{k: \bTheta^*_{jk}\neq 0\}. 
\eel
A simple version of our main theoretical result can be stated as follows. 

{\theorem\label{main} Let $\hbTheta$ and $\hbOmega$ be the scaled Lasso estimators 
defined in (\ref{joint-min}), 
(\ref{tTheta}) and (\ref{hTheta}) below with penalty level $\lam_0=A\sqrt{4(\log p)/n}$, $A>1$, 
based on $n$ iid observations from $N(0,\bSigma^*)$. 
Suppose the spectrum norm of $\bOmega^* = (\diag\bSigma^*)^{1/2}\bTheta^*(\diag\bSigma^*)^{1/2}$ 
is bounded and that $d^2(\log p)/n \to 0$. Then, 
\bes
\|\hbOmega-\bOmega^*\|_2 = O_P(1)d\sqrt{(\log p)/n}=o(1), 
\ees
where $\|\cdot\|_2$ is the spectrum norm (the $\ell_2$ matrix operator norm). 
If in addition the diagonal elements of $\bTheta^*$ is uniformly bounded, then 
\bes
\|\hbTheta-\bTheta^*\|_2 = O_P(1)d\sqrt{(\log p)/n}=o(1). 
\ees
}

Theorem \ref{main} provides a simple boundedness condition on the spectrum norm of 
$\bOmega^*$ for the convergence of $\hbOmega$ in spectrum norm with sample size $n\gg d^2\log p$. 
The additional condition on the diagonal of $\bTheta^*$ is natural due to scale change. 
The boundedness condition on the spectrum norm of 
$(\diag\bSigma^*)^{1/2}\bTheta^*(\diag\bSigma^*)^{1/2}$ and the diagonal of $\bTheta^*$ 
is weaker than the boundedness of the $\ell_1$ operator norm assumed in \cite{Yuan10,CaiLL11} 
since the boundedness of $\diag\bSigma^*$ is also needed there. 
When the ratio of the $\ell_1$ operator norm and spectrum norm of the precision matrix diverges to infinity, 
the proposed estimator has a faster proven convergence rate. 
This sharper result is a direct consequence of the faster convergence rate of the scaled Lasso estimator 
of the noise level in linear regression. 
To the best of our knowledge, it is unclear if the $\ell_1$ regularization method of \cite{Yuan10,CaiLL11} 
also achieve the convergence rate under the weaker spectrum norm condition.  

An important advantage of the scaled Lasso is that the penalty level is automatically 
set to achieve the optimal convergence rate in the regression model for the estimation of each column 
of the inverse matrix. 
This raises the possibility for the scaled Lasso to outperform methods using a single unscaled 
penalty level for the estimation of all columns such as the GLasso and CLIME. 
We provide an example in Section 7 to demonstrate the feasibility of such a scenario. 

Another contribution of this paper is to study the scaled Lasso at a smaller penalty 
level than those based on $\ell_\infty$ bounds of the noise. The $\ell_\infty$-based analysis requires 
a penalty level $\lam_0$ satisfying  $P\{N(0,1/n) > \lam_0/A\} = \eps/p$ for a small $\eps$ and $A>1$. 
For $A\approx 1$ and $\eps = p^{o(1)}$, this penalty level is comparable to the universal penalty level 
$\sqrt{(2/n)\log p}$. However, $\eps = o(1/p)$, or equivalently $\lam_0\approx \sqrt{(4/n)\log p}$, 
is required if the union bound is used to simultaneously control the error of $p$ applications of the scaled 
Lasso in the estimation of individual columns of a precision matrix. 
This may create a significant gap between theory and implementation. 
We close this gap by providing a theory based on a sparse $\ell_2$ measure of the noise, 
corresponding to a penalty level satisfying  $P\{N(0,1/n) > \lam_0/A\} = k/p$ with $A>1$ 
and a potentially large $k$. 
This penalty level provides a faster convergence rate than the universal penalty level 
in linear regression when $\log(p/k)\approx \log(p/\|\bbeta\|_0)\ll \log p$. 
Moreover, the new analysis provides a higher concentration of the error 
so that the same penalty level $\lam_0 \approx \sqrt{(2/n)\log(p/k)}$ can be used to simultaneously 
control the estimation error in $p$-applications of the scaled Lasso for the estimation of a precision matrix. 

The rest of the paper is organized as follows. 
In Section 2, we present the scaled Lasso method for the estimation of the inversion of a 
nonnegative definite matrix. 
In Section 3, we study the estimation error of the proposed method. 
In Section 4, we provide a theory for the Lasso and its scaled version with higher proven 
concentration at a smaller, practical penalty level. 
In Section 5, we study the least square estimation after the scaled Lasso selection. 
Simulation studies are presented in Section 6. 
In Section 7, we discuss the benefits of using the scaled penalty levels for the 
estimation of different columns of the precision matrix, compared with an optimal fixed penalty level 
for all columns. 
Section 8 provides all the proofs.

We use the following notation throughout the paper.
For real $x$, $x_+=\max(x,0)$.
For a vector $\bv=(v_1,\dots,v_p)$, $\|\bv\|_q=(\sum_j|v_j|^q)^{1/q}$ is the $\ell_q$ norm with the special $\|\bv\|=\|\bv\|_2$ and the usual extensions $\|\bv\|_{\infty}=\max_j |v_j|$ and $\|\bv\|_0=\#\{j:v_j\neq 0\}$.
For matrices $\bM$, $\bM_{i,*}$ is the $i$-th row and $\bM_{*,j}$ the $j$-th column,
$\bM_{A,B}$ represents the submatrix of $\bM$ with rows in $A$ and columns in $B$,
$\|\bM\|_q=\sup_{\|\bv\|_q=1}\|\bM\bv\|_q$ is the $\ell_q$ matrix operator norm.
In particular, $\|\cdot\|_2$ is the spectrum norm for symmetric matrices.
Moreover, we may denote the set $\{j\}$ by $j$ and denote the set $\{1,\dots, p\}\setminus \{j\}$ by $-j$ in the subscript.

\section{{Matrix inversion via scaled Lasso}}
Let $\bSigmabar$ be a nonnegative-definite data matrix and $\bTheta^*$ be a positive-definite target matrix 
with $\bSigmabar\bTheta^*\approx \bI$. 
In this section, we describe the relationship between positive-definite matrix inversion and linear regression 
and propose an estimator for $\bTheta^*$ via scaled Lasso, 
a joint convex minimization for the estimation of regression coefficients and noise level.

We use the scaled Lasso to estimate $\bTheta^*$ column by column. 
Define $\sigma_j>0$ and $\bbeta\in \real^{p\times p}$ by
\bel{sigma-beta-j}
\sigma_j^{2}=(\bTheta^*_{jj})^{-1},\quad
\bbeta_{*,j}=-\bTheta^*_{*,j}\sigma_j^{2}=-\bTheta^*_{*,j}(\bTheta^*_{jj})^{-1}.
\eel
In the matrix form, we have the following relationship
\bel{Theta}
\diag\bTheta^*=\diag(\sigma_{j}^{-2}, j=1,\dots,p), \quad\bTheta^*=-\bbeta(\diag\bTheta^*).
\eel
Let $\bSigma^*=(\bTheta^*)^{-1}$. 
Since $(\pa/\pa\bb_{-j})\bb^T\bSigma^*\bb = 2\bSigma^*_{-j,*}\bb=0$ at $\bb=\bbeta_{*,j}$,
one may estimate the $j$-th column of $\bbeta$ by minimizing the $\ell_1$ penalized quadratic loss. 
In order to penalize 
the unknown coefficients in the same scale, we adjust the $\ell_1$ penalty with 
diagonal standardization, leading to the following penalized quadratic loss:
\bel{Lasso}
\bb^T\bSigmabar\bb/2+\lam\sum_{k=1}^p \Sigmabar_{kk}^{1/2}|b_k|. 
\eel
For $\bSigmabar = \bX^T\bX/n$ and $b_j=-1$, $\bb^T\bSigmabar\bb = \|\bx_j - \sum_{k\neq j}b_k\bx_k\|_2^2/n$, 
so that (\ref{Lasso}) is the penalized loss for the Lasso in linear regression of $\bx_j$ against 
$\{\bx_k, k\neq j\}$. 
This is similar to the procedures in \cite{Yuan10,CaiLL11} that use the Dantzig selector 
to estimate $\bTheta^*_{*,j}$ column-by-column. 
However, one still needs to choose a penalty level $\lam$ and to estimate $\sigma_j$ 
in order to recover $\bTheta^*$ via (\ref{Theta}). 
A solution to resolve these two issues is the scaled Lasso \cite{SunZ11}:
\bel{joint-min}
\{\hbbeta_{*,j}, \hsigma_j\}
= \argmin_{\bb,\sigma}\Big\{\frac{\bb^T\bSigmabar\bb}{2\sigma} + \frac{\sigma}{2}
+\lam_0\sum_{k=1}^p \Sigmabar_{kk}^{1/2}|b_k|: b_j=-1\Big\}
\eel
with $\lam_0\approx \sqrt{(2/n)\log p}$.
The scaled Lasso (\ref{joint-min}) is a solution of joint convex minimization in $\{\bb,\sigma\}$ \cite{HuberR09,Antoniadis10}.
Since $\bbeta^T\bSigma^*\bbeta = (\diag\bTheta^*)^{-1}\bTheta^*(\diag\bTheta^*)^{-1}$,
\bes
\diag\big(\bbeta^T\bSigma^*\bbeta\big) = (\diag\bTheta^*)^{-1}
=\diag(\sigma^2_j, j=1,\ldots,p).
\ees
Thus, (\ref{joint-min}) is expected to yield consistent estimates of $\sigma_j=(\bTheta^*_{jj})^{-1/2}$.

An iterative algorithm has been provided by \cite{SunZ11} to compute the scaled Lasso estimator (\ref{joint-min}). We rewrite the algorithm in the form of matrices.
For each $j\in \{1,\dots,p\}$, the Lasso path is given by the estimates $\hbbeta_{-j,j}(\lam)$ satisfying the following Karush-Kuhn-Tucker conditions: for all $k\neq j$,
\bel{KKT}
\begin{cases}
\Sigmabar_{kk}^{-1/2}\bSigmabar_{k,*}\hbbeta_{*,j}(\lam)= -\lam\sgn(\hbeta_{k,j}(\lam)),& \hbeta_{k,j}\neq 0,
\cr \Sigmabar_{kk}^{-1/2}\bSigmabar_{k,*}\hbbeta_{*,j}(\lam) \in \lam[-1,1], & \hbeta_{k,j}=0,
\end{cases}
\eel
where $\hbbeta_{jj}(\lam)=-1$. Based on the Lasso path $\hbbeta_{*,j}(\lam)$, the scaled Lasso estimator $\{\hbbeta_{*,j}, \hsigma_j\}$ is computed iteratively by
\bel{alg}
\hsigma_j^2 \leftarrow \hbbeta_{*,j}^T\bSigmabar\hbbeta_{*,j}, \quad
\lam \leftarrow \hsigma_j\lam_0,\quad
\hbbeta_{*,j} \leftarrow \hbbeta_{*,j}(\lam).
\eel
Here the penalty level of the Lasso is determined by the data without using cross-validation.
We then simply take advantage of the relationship (\ref{Theta}) and compute the coefficients and noise levels by the scaled Lasso for each column
\bel{tTheta}
\diag\tbTheta=\diag(\hsigma_{j}^{-2}, j=1,\dots,p), \quad\tbTheta=-\hbbeta({\diag\tbTheta}).
\eel

Now we have constructed an estimator for $\bTheta^*$. In our primary example of taking $\bSigmabar$ 
as a sample covariance matrix, the target $\bTheta^*$ is the inverse covariance matrix. 
One may also be interested in estimating the inverse correlation matrix
\bel{Omega^*}
\bOmega^*=(\bR^*)^{-1}=\big\{\bD^{-1/2}\bSigma^*\bD^{-1/2}\big\}^{-1}
=\bD^{1/2}(\bSigma^*)^{-1}\bD^{1/2},
\eel
where $\bD=\diag(\bSigma^*)$ and $\bR^*=\bD^{-1/2}\bSigma^*\bD^{-1/2}$ is the population correlation matrix.
The diagonal matrix $\bD$ can be approximated by the diagonal of $\bSigmabar$.
Thus, the inverse correlation matrix is estimated by
\bel{tOmega}
\tbOmega&=&\hbD^{1/2}\tbTheta\hbD^{1/2} \text{ with } \hbD=\diag(\bSigmabar_{jj}, j=1,\dots,p).
\eel
The estimator $\tbOmega$ here is a result of normalizing the precision matrix estimator by the population variances.
Alternatively, we may estimate the inverse correlation matrix by using the population correlation matrix
\bes
\overline{\bR}=(\diag\bSigmabar)^{-1/2}\bSigmabar(\diag\bSigmabar)^{-1/2}
=\hbD^{-1/2}\bSigmabar\hbD^{-1/2}
\ees
as data matrix. Let $\{\hbalpha_{*,j},\htau_j\}$ be the solution of (\ref{joint-min}) with $\overline{\bR}$ in place of $\bSigmabar$. We combine these column estimators as in (\ref{tTheta}) to have an alternative estimator for $\bOmega^*$ as follows: 
\bes
\diag\big(\tbOmega^{\text{Alt}}\big)=\diag(\htau_{j}^{-2}, j=1,\dots,p), \quad 
\tbOmega^{\text{Alt}}=-\hbalpha\,\diag\big(\tbOmega^{\text{Alt}}\big).
\ees
Since $\overline{\bR}_{jj}=1$ for all $j$, it follows from (\ref{joint-min}) that
\bes
\hbalpha_{*,j}=\hbD^{1/2}\hbbeta_{*,j}\hbD^{-1/2}_{jj},\quad
\htau_j=\hsigma_j\hbD^{-1/2}_{jj}.
\ees
This implies
\bes
\tbOmega^{\text{Alt}}
=-\hbD^{1/2}\hbbeta{\diag(\hbD^{-1/2}_{jj}\hsigma_j^{-2}\hbD_{jj}, j=1,\dots,p)}
=\hbD^{1/2}\tbTheta\hbD^{1/2}=\tbOmega.
\ees
Thus, in this scaled Lasso approach, the estimator based on the normalized data matrix 
is exactly the same as the one based on the original data matrix followed by a normalization step. 
The scaled Lasso methodology is scale-free in the noise level, and as a result, the estimator for inverse 
correlation matrix is also scale free in diagonal normalization.

It is noticed that a good estimator for $\bTheta^*$ or $\bOmega^*$ should be a symmetric matrix. However, the estimators $\tbTheta$ and $\tbOmega$ do not have to be symmetric. We improve them by using a symmetrization step as in \cite{Yuan10},
\bel{hTheta}
\hbTheta=\argmin_{\bM: \bM^T=\bM}\|\bM-\tbTheta\|_1,\quad
\hbOmega=\argmin_{\bM: \bM^T=\bM}\|\bM-\tbOmega\|_1,
\eel
which can be solved by linear programming. It is obvious that $\hbTheta$ and $\hbOmega$ are both symmetric, 
but not guaranteed to be positive-definite. It follows from Theorem \ref{main} that 
$\hbTheta$ and $\hbOmega$ are positive-definite with large probability. 
Alternatively, semidefinite programming, which is somewhat more expensive computationally, 
can be used to produce a nonnegative-definite $\hbTheta$ and $\hbOmega$ in (\ref{hTheta}).

According to the definition, the estimators $\hbTheta$ and $\hbOmega$ have the same $\ell_1$ 
error rate as $\tbTheta$ and $\tbOmega$ respectively. 
A nice property of symmetric matrices is that the spectrum norm is bounded by the $\ell_1$ matrix norm. 
The $\ell_1$ matrix norm can be expressed more explicitly as the maximum $\ell_1$ norm of the 
columns, while the $\ell_\infty$ matrix norm is the maximum $\ell_1$ norm of the rows. 
Hence, for any symmetric matrix, the $\ell_1$ matrix norm is equivalent to the $\ell_\infty$ matrix norm, 
and the spectrum norm can be bounded by either of them. 
Since our estimators and target matrices are all 
symmetric, the error bound based on the spectrum norm could be studied by bounding the $\ell_1$ error 
as typically done in the existing literature. We will study the estimation error of (\ref{hTheta}) in Section~3.

To sum up, we propose to estimate the matrix inversion by (\ref{joint-min}), 
(\ref{tTheta}) and (\ref{hTheta}). 
The iterative algorithm (\ref{alg}) computes (\ref{joint-min}) 
based on a Lasso path determined by (\ref{KKT}). 
Then (\ref{tTheta}) translates the resulting estimators of (\ref{alg}) to column estimators 
and thus a preliminary matrix estimator is constructed. 
Finally, the symmetrization step (\ref{hTheta}) produces a symmetric estimate for our target matrix.

\section{Theoretical properties}
From now on, we suppose that the data matrix is the sample covariance matrix $\bSigmabar=\bX^T\bX/n$, 
where the rows of $\bX$ are iid $N(0,\bSigma^*)$. 
Let $\bTheta^* = (\bSigma^*)^{-1}$ be the precision matrix as the inverse of the population covariance 
matrix. Let $\bD$ be the diagonal of $\bSigma^*$, $\bR^* = \bD^{-1/2}\bSigma^*\bD^{-1/2}$ 
the population correlation matrix, $\bOmega^*=(\bR^*)^{-1}$ its inverse as in (\ref{Omega^*}). 
In this section, we study $\hbOmega$ and $\hbTheta$ in (\ref{hTheta}), respectively 
for the estimation of $\bOmega^*$ and $\bTheta^*$.

We consider a certain capped $\ell_1$ sparsity for individual columns of the inverse matrix as follows. 
For a certain $\eps_0>0$, a threshold level $\lam_{*,0}>0$ not depending on $j$ and an index set 
$S_j\subset \{1,\ldots,p\}\setminus \{j\}$, the capped $\ell_1$ sparsity condition measures the 
complexity of the $j$-th column of $\bOmega^*$ by 
\bel{s-starj}
|S_j| + (1-\eps_0)^{-1}\sum_{k\neq j, {k \not\in S_j} }
\frac{|\Omega_{kj}^*|}{(\Omega_{jj}^*)^{1/2}\lam_{*,0}} \le s_{*,j}. 
\eel
The condition can be written as 
\bes
\sum_{j\neq k}\min\left\{\frac{|\Omega_{kj}^*|}{(1-\eps_0)(\Omega_{jj}^*)^{1/2}\lam_{*,0}},1\right\}\le s_{*,j}
\ees 
if we do not care about the choice of $S_j$.  
In the $\ell_0$ sparsity case of $S_j=\{k: k\neq j, \Omega_{kj}^*\neq 0\}$, 
we may set $s_{*,j}=|S_j|+1$ as the degree for the $j$-th node in the graph induced by 
matrix $\Omega^*$ (or $\bTheta^*$). 
In this case, $d = \max_j (1+|S_j|)$ is the maximum degree as in (\ref{degree}). 

In addition to the sparsity condition on the inverse matrix, we also require a certain 
invertibility condition on $\bR^*$. Let $S_j\subseteq B_j\subseteq \{1,\ldots,p\}\setminus \{j\}$. 
A simple version of the required invertibility condition can be written as 
\bel{partial-inv}
\inf\left\{ \frac{\bu^T(\bR^*_{-j,-j})\bu}{\|\bu_{B_j}\|_2^2}: \bu_{B_j}\neq 0\right\} \ge c_*
\eel
with a fixed constant $c_*>0$. 
This condition requires a certain partial invertibility of the population correlation matrix. 
It certainly holds if the smallest eigenvalue of $\bR^*_{-j,-j}$ is no smaller than $c_*$ for all $j\le p$,  
or the spectrum norm of $\bOmega^*$ is no greater than $1/c_*$ as assumed in Theorem \ref{main}. 
In the proof of Theorems \ref{main-1} and \ref{thm-corrmatrix}, we only use a weaker 
version of condition (\ref{partial-inv}) in the form of (\ref{c_*}) with $\{\bSigma^*,\bSigmabar\}$ replaced by 
$\{\bR^*_{-j,-j},{\overline \bR}_{-j,-j}\}$ there.

{\theorem\label{main-1} Suppose $\bSigmabar$ is the sample covariance matrix of 
$n$ iid $N(0,\bSigma^*)$ vectors. 
Let $\bTheta^*=(\bSigma^*)^{-1}$ and $\bOmega^*$ as in (\ref{Omega^*}) be the inverses of 
the population covariance and correlation matrices. 
Let $\hbTheta$ and $\hbOmega$ be their scaled Lasso estimators defined in (\ref{joint-min}), 
(\ref{tTheta}) and (\ref{hTheta}) with a penalty level $\lam_0=A\sqrt{4(\log p)/n}$, $A>1$. 
Suppose (\ref{s-starj}) and (\ref{partial-inv}) hold 
with $\eps_0=0$ and $\max_{j\le p}(1+s_{*,j})\lam_0\le c_0$ for a certain constant $c_0>0$ depending on $c_*$ 
only. Then, the spectrum norm of the errors are bounded by 
\bel{bound-hbTheta}
&& \|\hbTheta-\bTheta^*\|_2 \le \|\hbTheta-\bTheta^*\|_1 
\le C\Big(\max_{j\le p}(\big\|\bD_{-j}^{-1}\|_\infty\Theta^*_{jj})^{1/2}s_{*,j}\lam_0 
+\big\|\bTheta^*\big\|_1\lam_0\Big),\qquad
\\ \label{bound-hbOmega} && \|\hbOmega-\bOmega^*\|_2 \le \|\hbOmega-\bOmega^*\|_1 
\le C\Big(\max_{j\le p} (\Omega_{jj}^*)^{1/2}s_{*,j}\lam_0+\|\bOmega^*\|_1\lam_0\Big), 
\eel
with large probability, where $C$ is a constant depending on $\{c_0,c_*,A\}$ only. 
Moreover, the term $\|\bTheta^*\|_1\lam_0$ in (\ref{bound-hbTheta}) can be replaced by 
\bel{main-1-3}
&& \max_{j\le p}\|\bTheta_{*,j}\|_1s_{*,j}\lam_0^2+\tau_n(\bTheta^*),
\eel
where $\tau_n(\bM) = \inf\{\tau: \sum_j \exp( - n\tau^2/\|\bM_{*,j}\|_1^2)\le 1/e\}$. 
}

\medskip
Theorem \ref{main-1} implies Theorem 1 due to $s_{*,j}\le d-1$, 
$1/D_{jj} \le \Theta^*_{jj}\le \|\bTheta^*\|_2$, $\|\bTheta^*\|_1\le d\max_j \Theta^*_{jj}$ 
and similar inequalities for $\bOmega^*$. We note that $B_j=S_j$ in (\ref{partial-inv}) gives the 
largest $c_*$ and thus the sharpest error bounds in Theorem \ref{main-1}. 
In Section 7, we give an example to demonstrate the advantage of this theorem. 

In a 2011 arXiv version of this paper (http://arxiv.org/pdf/1202.2723v1.pdf), 
we are able to demonstrate good numerical performance of the scaled Lasso estimator 
with the universal penalty level $\lam_0 = \sqrt{2(\log p)/n}$, compared with some 
existing methods, but not the larger penalty level $\lam_0 > \sqrt{4(\log p)/n}$ in Theorems \ref{main} and \ref{main-1}. 
Since a main advantage of our proposal is automatic selection of the penalty level without 
resorting to cross validation, a question arises as to whether a theory can be developed 
for a smaller penalty level to match the choice in a demonstration of good performance 
of the scaled Lasso in our simulation experiments. 

We are able to provide an affirmative answer in this version of the paper 
by proving a higher concentration of the error of the scaled Lasso at a smaller penalty 
level as follows. 
Let $L_n(t)$ be the $N(0,1/n)$ quantile function satisfying 
\bes
P\big\{N(0,1) > n^{1/2}L_n(t)\big\}=t.
\ees
Our earlier analysis is based on existing oracle inequalities of the Lasso which 
holds with probability $1-2\eps$ when the inner product of design vectors 
and noise are bounded by their $\eps/p$ and $1-\eps/p$ quantiles. 
Application of the union bound in $p$-application of the Lasso requires a threshold level 
$\lam_{*,0}=L_n(\eps/p^2)$ with a small $\eps>0$, which matches $\sqrt{4(\log p)/n}$ 
with $\eps\asymp 1/\sqrt{\log p}$ in Theorems \ref{main} and \ref{main-1}. 
Our new analysis of the scaled Lasso allows a threshold level 
\bes
\lam_{*,0} = L_{n-3/2}(k/p)
\ees
with $k\asymp s\log(p/s)$, where $s = 1+\max_j s_{*,j}$. More precisely, 
we require a penalty level $\lam_0\ge A\lam_{*,0}$ with a constant $A$ satisfying 
\bel{cond-xi1}
A-1>A_1\ge \max_j \Big\{ \Big[\frac{e^{1/(4n-6)^2}4k}{m_j(L^4+2L^2)}\Big]^{1/2}
+ \frac{e^{1/(4n-6)^2}}{L\sqrt{2\pi}}\sqrt{\psi_j}
+ \frac{L_1(\eps/p)}{L}\sqrt{\psi_j}\Big\}, 
\eel
where $L=L_1(k/p)$, $s_{*,j}\le m_j\le \min(|B_j|,C_0 s_{*,j})$ with the $s_{*,j}$ and $B_j$ 
in (\ref{s-starj}) and (\ref{partial-inv}), and 
$\psi_j=\kappa_+(m_j;{\bR}_{-j,-j})/m_j + L_n(5\eps/p^2)$
with 
\bel{sparse-eigen}
\kappa_+(m;\bSigma) = \max_{\|\bu\|_0=m,\|\bu\|_2=1}\bu^T\bSigma\bu.
\eel

{\theorem\label{thm-corrmatrix} 
Let $\{\bSigmabar, \bSigma^*, \bTheta^*,\bOmega^*\}$ be matrices as in Theorem \ref{main-1}, 
and $\hbTheta$ and $\hbOmega$ be the scaled Lasso estimators 
with a penalty level $\lam_0\ge A\lam_{*,0}$ where $\lam_{*,0} = L_{n-3/2}(k/p)$. 
Suppose (\ref{s-starj}) and (\ref{partial-inv}) hold with certain $\{S_j,s_{*,j},\eps_0,B_j,c_*\}$, 
(\ref{cond-xi1}) holds with constants $\{A, A_1,C_0\}$ and certain integers $m_j$, and 
$P\{(1-\eps_0)^2\le \chi_n^2/n\le (1+\eps_0)^2 \}\le \eps/p$. Then,  
there exist constants $c_0$ depending on $c_*$ only and 
$C$ depending on $\{A, A_1, C_0, c_*,c_0\}$ only such that 
when $\max_j s_{*,j} \lam_0 \le c_0$, 
the conclusions of Theorem \ref{main-1} hold 
with at least probability $1-6\eps-2k\sum_j(p-1-|B_j|)/p$. 
}
\medskip

The condition $\max_j s_{*,j} \lam_0 \le c_0$ on (\ref{s-starj}), which controls the 
capped $\ell_1$ sparsity of the inverse correlation matrix, weakens the $\ell_0$ sparsity condition 
$d\sqrt{(\log p)/n}\to 0$. 

The extra condition on the upper sparse eigenvalue $\kappa_+(m;\bR^*_{-j,-j})$ in (\ref{cond-xi1}) is mild, 
since it only requires a small $\kappa_+(m;\bR^*)/m$ that is actually decreasing in $m$.

The invertibility condition (\ref{partial-inv}) is used to regularize the design matrix in linear regression 
procedures. As we mentioned earlier, condition (\ref{partial-inv}) holds if 
the spectrum norm of $\bOmega^*$ is bounded by $1/c_*$. 
Since $(\bR^*)^{-1}=\bOmega^*=(\diag\bSigma^*)^{1/2}\bTheta^*(\diag\bSigma^*)^{1/2}$, 
it suffices to have 
\bes
\|(\bR^*)^{-1}\|_2 \le 
\max\bSigma^*_{jj}\|\bTheta^*\|_2\le 1/c_*.
\ees
To achieve the convergence rate $d\sqrt{(\log p)/n}$, 
both \cite{Yuan10,CaiLL11} require conditions $\|\bTheta^*\|_1=O(1)$ and $\max\bSigma^*_{jj} =O(1)$. 
In comparison, the spectrum norm condition is not only weaker than the $\ell_1$ operator norm 
condition, but also more natural for the convergence in spectrum norm.

Our sharper theoretical results are consequences of using the scaled Lasso 
estimator (\ref{joint-min}) and its fast convergence rate in linear regression. 
In \cite{SunZ11}, a convergence rate of order $s_*(\log p)/n$ was established for the scaled Lasso estimation 
of the noise level, compared with an oracle noise level as the moment estimator based on the noise vector. 
In the context of the column-by-column application of the scaled Lasso for precision matrix estimation, 
the results in \cite{SunZ11} can be written as  
\bel{oracle-scale}
\Big|\frac{\sigma_j^*}{\hsigma_j}-1\Big| \le C_1 s_{*,j}\lam_{0}^2,\quad
\sum_{k\neq j}\Sigmabar_{kk}^{1/2}|\hbeta_{k,j}-\beta_{k,j}|\sqrt{\Theta^*_{jj}} \le C_2s_{*,j}\lam_{0}, 
\eel
where $\sigma_j^* = \|\bX\bbeta_{*,j}\|_2/\sqrt{n}$. We note that 
$n(\sigma_j^*)^2\Theta^*_{jj}$ is a chi-square variable with $n$ degrees of freedom 
when $\bX$ has iid $N(0,\bSigma^*)$ rows. 
The oracle inequalities in (\ref{oracle-scale}) play a crucial role in our analysis of 
the proposed estimators for inverse matrices, as the following proposition attests. 

\begin{proposition}\label{prop-matrix}
Let $\bTheta^*$ be a nonnegative definite target matrix, $\bSigma^*=(\bTheta^*)^{-1}$, 
and $\bbeta=-\bTheta^*(\diag\bTheta^*)^{-1}$. 
Let $\hbTheta$ and $\hbOmega$ be defined as (\ref{tTheta}) and (\ref{hTheta}) 
based on certain $\hbbeta$ and $\hsigma_j$ satisfying (\ref{oracle-scale}). 
Suppose further that 
\bel{prop-matrix-1}
|\Theta_{jj}^*(\sigma^*_j)^2-1|\le C_0\lam_0,\  
\max_j|(\Sigmabar_{jj}/\Sigma^*_{jj})^{-1/2}-1|\le C_0\lam_0, 
\eel
and that 
$\max\{4C_0\lam_0,4\lam_0,C_1s_{*,j}\lam_0\}\le 1$. 
Then, (\ref{bound-hbTheta}) and (\ref{bound-hbOmega}) hold with a constant 
$C$ depending on $\{C_0,C_2\}$ only. Moreover, if $n\bTheta_{jj}^*(\sigma^*_j)^2\sim \chi^2_n$, 
then the term $\lam_0\|\bTheta^*\|_1$ in (\ref{bound-hbTheta}) can be replaced by (\ref{main-1-3}) 
with large probability. 
\end{proposition}

While the results in \cite{SunZ11} requires a penalty level $A\sqrt{(2/n)\log(p^2)}$  
to allow simultaneous application of (\ref{oracle-scale}) for all $j\le p$ via the union bound in proving  
Theorem~\ref{main-1},   
Theorem \ref{thm-corrmatrix} allows a smaller penalty level $\lam_{*,0}=A L_{n-3/2}(k/p)$ with $A>1$ and 
a potentially large $k\asymp s\log(p/s)$. 
This is based on new theoretical results for the Lasso and scaled Lasso 
developed in Section 4.

\section{Linear regression revisited}
This section provides certain new error bounds for the Lasso and scaled Lasso in the linear regression model.
Compared with existing error bounds, the new results characterize the concentration of the estimation and
prediction errors at fixed, smaller threshold levels. 
The new results also allow high correlation among certain nuisance design vectors.

Consider the linear regression model with standardized design and normal error:
\bes
\by = \bX\bbeta + \bep,\ \|\bx_j\|_2^2=n,\ \bep\sim N(0,\sigma^2\bI_n).
\ees
Let $\lam_{univ}=\sqrt{(2/n)\log p}$ be the universal penalty level \cite{DonohoJ94}.
For the estimation of $\bbeta$ and variable selection,
existing theoretical results with $p\gg n$ typically require a penalty level $\lam = A\sigma\lam_{univ}$,
with $A>1$, to guarantee rate optimality of regularized estimators.
This includes the scaled Lasso with a jointly estimated $\sigma$.
For the Dantzig selector \cite{CandesT07}, performance bounds have been established for $A=1$.

It is well understood that $\sigma\lam_{univ}$ in such theorems is
a convenient probabilistic upper bound of $\|\bX^T\bep/n\|_\infty$ for controlling the
maximum gradient of the squared loss $\|\by-\bX\bb\|_2^2/(2n)$ at $\bb=\hbbeta$.
For $\lam < \|\bX^T\bep/n\|_\infty$, variable selection is known to be inconsistent for the Lasso
and most other regularized estimates of $\bbeta$, and the analysis of such procedures
become more complicated due to false selection.
However, this does not preclude the possibility that such a smaller $\lam$ outperforms the
theoretical $\lam \ge \sigma\lam_{univ}$ for the estimation of $\bbeta$ or prediction. 

In addition to theoretical studies, a large volume of numerical comparisons among
regularized estimators exists in the literature.
In such numerical studies, the choice of penalty level is typically delegated to
computationally more expensive cross-validation methods.
Since cross-validation aims to optimize prediction performance, it may lead to a smaller penalty
level than $\lam=\sigma\lam_{univ}$. However, this gap between
$\lam \ge \sigma\lam_{univ}$ in theoretical studies and the possible choice of
$\lam <\sigma\lam_{univ}$ in numerical studies is largely ignored in the existing literature.

The purpose of this section is to provide rate optimal oracle inequalities for the Lasso and its
scaled version, which hold with at least probability $1-\eps/p$ for a reasonably small $\eps$,
at a fixed penalty level $\lam$ satisfying $P\{N(0,\sigma^2/n) > \lam/A\} = k/p$, with a given $A > 1$ 
and potentially large $k$, up to $k/(2\log(p/k))^2 \asymp s_*$, 
where $s_*$ is a complexity measure of $\bbeta$, e.g. $s_*=\|\bbeta\|_0$.

When the (scaled) Lasso is simultaneously applied to $p$ subproblems as in the case of matrix
estimation, the new oracle inequalities allow the use of the union bound to uniformly control the
estimation error in subproblems at the same penalty level.

Rate optimal oracle inequalities have been established for $\ell_1$ and concave regularized estimators 
in \cite{Zhang10,YeZ10} for penalty level $\lam=A \sigma \sqrt{c^*(2/n)\log(p/(\eps s_*))}$, 
where $c^*$ is an upper sparse eigenvalue, $A>1$ and $1-\eps$ is the guaranteed probability for 
the oracle inequality to hold. 
The new oracle inequalities remove the factors $c^*$ and $\eps$ from the penalty level, 
as long as $1/\eps$ is polynomial in $p$. 
The penalty level $A\sigma\sqrt{(2/n)\log(p/(\eps s))}$ has been considered for models of size
$s$ under $\ell_0$ regularization \cite{BirgeM01,BirgeM07,BuneaTW07-Aggr, AbramovichG10}. 

To bound the effect of the noise when $\lam < \|\bX^T\bep/n\|_\infty$, 
we use a certain sparse $\ell_2$ 
norm to control the excess of $\bX^T\bep/n$ over a threshold level $\lam_*$. 
The sparse $\ell_q$ norm was used in the analysis of regularized estimators before 
\cite{CandesT07,ZhangH08,Zhang09-l1,Zhang10,CaiWX10,YeZ10}, 
but it was done without a formal definition of the quantity to the best of our knowledge. 
To avoid repeating existing calculation, 
we define the norm and its dual here and summarize their properties in a proposition.

For $1\le q\le \infty$ and $t>0$, 
the sparse $\ell_q$ norm and its dual are defined as
\bel{sparse-norm}
&& \|\bv\|_{(q,t)} = \max_{|B| < t+1}\|\bv_B\|_q,\ \|\bv\|_{(q,t)}^* = \max_{\|\bu\|_{(q,t)}\le 1}\bu^T\bv. 
\eel
The following proposition describes some of their basic properties. 

\begin{proposition}\label{prop-sparse-norm} 
Let $m\ge 1$ be an integer, $q'=q/(q-1)$ and $a_q=(1-1/q)/q^{1/(q-1)}$. \\
(i) Properties of $\|\cdot\|_{(q,{m})}$: $\|\bv\|_{(q,{m})}\downarrow q$, $\|\bv\|_{(q,{m})}/{m}^{1/q} \downarrow {m}$,
$\|\bv\|_{(q,{m})}/{m}^{1/q} \uparrow q$,
\bel{prop-sparse-norm-1}
\|\bv\|_\infty = \|\bv\|_{(q,1)} \le \|\bv\|_{(q,{m})} \le (\|\bv\|_q)\wedge({m}^{1/q}\|\bv\|_\infty), 
\eel
and $\|\bv\|_q^q\le \|\bv\|_{(q,{m})}^q+(a_q/{m})^{q-1}\|\bv\|_1^q$. \\
(ii) Properties of $\|\cdot\|_{(q,{m})}^*$: $m^{1/q}\|\bv\|_{(q,{m})}^*\downarrow q$, and 
\bel{prop-sparse-norm-2}
\max\big(\|\bv\|_{q'},{m}^{-1/q}\|\bv\|_1\big) \le \|\bv\|_{(q,{m})}^* \le 
\min\big(\|\bv\|_{(q',{m}/a_q)} + {m}^{-1/q}\|\bv\|_1,\|\bv\|_1\big). 
\eel
(iii) Let $\bSigmabar = \bX^T\bX/n$ and $\kappa_+({m};\bM)$ be the sparse eigenvalue in 
(\ref{sparse-eigen}). Then,
\bes
\|\bSigmabar\bv\|_{(2,{m})} 
\le \min\Big\{\kappa_+^{1/2}({m};\bSigmabar)\|\bSigmabar^{1/2}\bv\|_{2}, 
\kappa_+({m};\bSigmabar)\|\bv\|_{2}\Big\}. 
\ees
\end{proposition}

\subsection{Lasso with smaller penalty: analytical bounds}
The Lasso path is defined as an $\R^p$-valued function of $\lam>0$ as
\bel{lasso-small-lam}
\hbbeta(\lam) = \argmin_{\bb}\Big\{\|\by-\bX\bb\|_2^2/(2n)+\lam\|\bb\|_1\Big\}.
\eel
For threshold levels $\lam_*>0$, we consider $\bbeta$ satisfying the following complexity bound,
\bel{s_*}
|S| + \sum_{j\not\in S}|\beta_j|/\lam_* \le s_*
\eel
with a certain $S\subset\{1,\ldots,p\}$. 
This includes the $\ell_0$ sparsity condition $\|\bbeta\|_0 = s_*$ with $S=\supp(\bbeta)$
and allows $\|\bbeta\|_0$ to far exceed $s_*$ with many small $|\beta_j|$.

The sparse $\ell_2$ norm of a soft-thresholded vector $\bv$, at threshold level $\lam_*$ in (\ref{s_*}), is 
\bel{zeta}
\zeta_{(2,m)}(\bv,\lam_*) = \|(|\bv|-\lam_*)_+\|_{(2,m)} = 
\max_{|J|\le m}\Big\{ \sum_{j\in J}(|v_j| - \lam_*)_+^2\Big\}^{1/2}. 
\eel
Let $B\subseteq \{1,\ldots,p\}$ and 
\bel{oeps}
\bz = (z_1,\ldots, z_p)^T = \bX^T\bep/n. 
\eel
We bound the effect of the excess of the noise over $\lam_*$ under the condition 
\bel{noise-cond}
\|\bz_{B^c}\|_\infty\le \lam_*,\ \zeta_{(2,m)}(\bz_B,\lam_*) \le A_1m^{1/2}\lam_*, 
\eel
 for some $A_1\ge 0$. We prove that when $\lam\ge A\lam_*$ with $A>1+A_1$ and (\ref{noise-cond}) holds, 
a scaled version of $\hbbeta(\lam)-\bbeta$ belongs to a set 
$\scrU(\bSigmabar,S,B;A,A_1,m,s_*-|S|)$, where 
\bel{U}
&& 
\scrU(\bSigma,S,B;A,A_1,m,m_1) 
\\ \nonumber &=& \Big\{\bu: \bu^T\bSigma\bu+(A-1)\|\bu_{S^c}\|_1
\le (A+1)\|\bu_S\|_1 + A_1m^{1/2}\|\bu_B\|_{(2,m)}^*  + 2A m_1\Big\}. 
\eel
This leads to the definition of 
\bel{M_pred}
M^*_{pred} = \sup\Big\{\frac{\bu^T\bSigmabar\bu/A^2}{m_1+|S|}: 
\bu\in \scrU(\bSigmabar,S,B;A,A_1,m,m_1)\Big\}
\eel
as a constant factor for the prediction error of the Lasso and 
\bel{M_q}
M^*_{q} = \sup\Big\{\frac{\|\bu\|_q/A}{(m_1+|S|)^{1/q}}: 
\bu\in \scrU(\bSigmabar,S,B;A,A_1,m,m_1) \Big\}
\eel
for the $\ell_q$ estimation error of the Lasso. 

The following theorem provides analytic error bounds for the Lasso prediction and estimation
under the sparse $\ell_2$ norm condition (\ref{noise-cond}) on the noise.
This is different from existing analyses of the Lasso based on the $\ell_\infty$ noise bound
$\|\bX^T\bep/n\|_\infty \le \lam_*$. 
In the case of Gaussian error, 
(\ref{noise-cond}) allows a fixed threshold level $\lam_* = \sigma\sqrt{(2/n)\log(p/m)}$ to uniformly
control the error of $p$ applications of the Lasso for the estimation of a precision matrix.
When $m\asymp s_*$ and $\sigma\sqrt{(2/n)\log(p/m)} \ll \sigma\sqrt{(2/n)\log p}$, 
using such smaller $\lam_*$ is necessary for achieving error bounds 
with the sharper rate corresponding to $\sigma\sqrt{(2/n)\log(p/m)}$.

\begin{theorem}\label{thm-small-lam}
Suppose (\ref{s_*}) holds with certain $\{S,s_*,\lam_*\}$.  
Let $A>1$, $\hbbeta=\hbbeta(\lam)$ be the Lasso estimator with penalty level $\lam\ge A\lam_*$, 
$\bh=\hbbeta-\bbeta$, and $m_1=s_*-|S|$. 
If (\ref{noise-cond}) holds with $A_1\ge 0$, a positive integer $m$ and $B\subseteq \{1,\ldots,p\}$, then
\bel{th-small-lam-pred-est}
\|\bX\bh\|_2^2/n \le M^*_{pred} s_* \lam^2,\quad 
\|\bh\|_q \le M^*_{q} s_*^{1/q}\lam. 
\eel
\end{theorem}

\begin{remark}\label{remark-cone} 
Theorem \ref{thm-small-lam} covers $\|\bX^T\bep/n\|_\infty\le \lam_*$ as a special case with $A_1=0$. 
In this case, the set (\ref{U}) does not depend on $\{m,B\}$. 
For $A_1=0$ and $|S|=s_*$ ($m_1=0$), (\ref{U}) contains all vectors satisfying 
a basic inequality $\bu^T\bSigma\bu+(A-1)\|\bu_{S^c}\|_1 \le (A+1)\|\bu_S\|_1$
\cite{BickelRT09,vandeGeerB09,YeZ10} and 
Theorem \ref{thm-small-lam} still holds when (\ref{U}) is replaced by the smaller 
\bel{S-cone}
\scrU_-(\bSigma,S,A) 
= \Big\{\bu: \|\bSigmabar_{S,*}\bu\|_\infty \le A+1, 
u_j\bSigmabar_{j,*}\bu\le -|u_j|(A-1)\ \forall j\not\in S\big\}. 
\eel
Thus, in what follows, we always treat $\scrU(\bSigma,S,B;A,0,m,0)$ as 
$\scrU_-(\bSigma,S,A)$ when $A_1=0$ and $|S|=s_*$. 
This yields smaller constants $\{M^*_{pred},M^*_q\}$ in (\ref{M_pred}) and (\ref{M_q}).  
\end{remark}

The purpose of including a choice $B$ in (\ref{noise-cond}) 
is to achieve bounded $\{M^*_{pred},M^*_1\}$ in the presence of some highly correlated 
design vectors outside $S\cup B$ when $\bSigmabar_{S\cup B,(S\cup B)^c}$ is small. 
Since $\|\bu_B\|_{(2,m)}^*$ is increasing in $B$, a larger $B$ leads to a larger set (\ref{U}) 
and larger $\{M^*_{pred},M^*_q\}$. 
However, (\ref{noise-cond}) with smaller $B$ typically requires larger $\lam_*$. 
Fortunately, the difference in the required $\lam_*$ in (\ref{noise-cond}) is of smaller order than 
$\lam_*$ between the largest $B=\{1,\ldots,p\}$ and smaller $B$ with $|B^c|\le p/m$. 
We discuss the relationship between $\{M^*_{pred},M^*_q\}$ and existing 
conditions on the design in the next section, 
along with some simple upper bounds for $\{M^*_{pred},M^*_1,M^*_2\}$. 

\subsection{Scaled Lasso with smaller penalty: analytical bounds}
The scaled Lasso estimator is defined as
\bel{scaledlasso}
\{\hbbeta,\hsigma\} = \argmin_{\bb,\sigma}\Big\{\|\by-\bX\bb\|_2^2/(2n\sigma)+\lam_0\|\bb\|_1+\sigma/2\Big\},
\eel
where $\lam_0>0$ is a scale-free penalty level. 
In this section, we describe the implication of Theorem \ref{thm-small-lam} on the scaled Lasso. 

A scaled version of (\ref{s_*}) is 
\bel{s_*0}
|S| + \sum_{j\not\in S}|\beta_j|/(\sigma^*\lam_{*,0}) = s_{*,0} \le s_*, 
\eel
where $\sigma^*=\|\bep\|_2/\sqrt{n}$ is an oracle estimate of the noise level and 
$\lam_{*,0}>0$ is a scaled threshold level. 
This holds automatically under (\ref{s_*}) when $S\supseteq \supp(\bbeta)$. 
When $\bbeta_{S^c}\neq 0$, (\ref{s_*0}) can be viewed as an event of large probability. 
When 
\bel{s_*0-bd}
|S|+(1-\eps_0)^{-1}\sum_{j\not\in S}\frac{|\beta_j|}{\sigma\lam_{*,0}}\le s_*
\eel
and $\bep\sim N(0,\sigma \bI_n)$, 
$P\big\{s_{*,0}\le s_*\big\} \ge P\big\{\chi^2_n/n \ge (1-\eps_0)^2\big\} \to 1$ 
for fixed $\eps_0>0$. Let 
\bel{M_sigma}
M^*_\sigma = \sup_{\bu\in\scrU}
\left\{ \frac{\bu^T\bSigmabar\bu}{s_*A^2}+\frac{2\|\bu\|_1}{s_*A^2}
+\frac{2A_1m^{1/2}\|\bu_B\|_{(2,m)}^*}{s_*A^2} \right\}
\eel
with $\scrU = \scrU(\bSigmabar,S,B;A,A_1,m,m_1)$ in (\ref{U}), as in (\ref{M_pred}) and (\ref{M_q}). Set 
\bes
\eta_* = M^*_\sigma A^2\lam_{*,0}^2s_*,\ 
\lam_0\ge A\lam_{*,0}/\sqrt{(1-\eta_*)_+},\ 
\eta_0 = M_\sigma^*\lam_0^2 s_*. 
\ees

\begin{theorem}\label{thm-scaled}
Suppose $\eta_0<1$. 
Let $\{\hbbeta,\hsigma\}$ be the scaled Lasso estimator in (\ref{scaledlasso}),
$\phi_1 = 1/\sqrt{1+\eta_0}$, $\phi_2 =1/\sqrt{1-\eta_0}$, 
and $\sigma^*=\|\bep\|_2/\sqrt{n}$. 
Suppose (\ref{noise-cond}) holds with 
$\{\bz_B,\lam_*\}$ replaced by $\{\bz_B/\sigma^*,\lam_{*,0}\}$. \\
(i) Let $\bh^* = (\hbbeta - \bbeta)/\sigma^*$. Suppose (\ref{s_*0}) holds. Then, 
\bel{thm-scaled-1}
\phi_1 < \hsigma/\sigma^* < \phi_2,\
\|\bX\bh^*\|_2^2/n < M_{pred}^*s_*(\phi_2\lam_0)^2,\
\|\bh^*\|_q < M_{q}^*s_*\phi_2\lam_0. 
\eel
(ii) Let $\bh = \hbbeta - \bbeta$. Suppose (\ref{s_*0-bd}) holds and 
$1-\eps_0 \le \sigma^*/\sigma\le1+\eps_0$. Then, 
\bel{thm-scaled-2}
&& (1-\eps_0)\phi_1 < \hsigma/\sigma < \phi_2(1+\eps_0),\
\\ \nonumber && \|\bX\bh\|_2^2/n < (1+\eps_0)^2M_{pred}^*s_*(\sigma\phi_2\lam_0)^2,\
\\ \nonumber &&\|\bh\|_q < (1+\eps_0)M_{q}^*s_*\sigma \phi_2\lam_0. 
\eel
\end{theorem}

Compared with Theorem \ref{thm-small-lam},
Theorem \ref{thm-scaled} requires nearly identical conditions on the design $\bX$, 
the noise and penalty level under proper scale.
It essentially allows the substitution of $\{\by,\bX,\bbeta\}$ by 
$\{\by/\sigma^*,\bX,\bbeta/\sigma^*\}$ when $\eta_0$ is small. 

Theorems \ref{thm-small-lam} and \ref{thm-scaled} require an upper bound (\ref{noise-cond}) 
for the sparse $\ell_2$ norm of the excess noise as well as upper bounds for the constant factors 
$\{M^*_{pred},M^*_q,M^*_\sigma\}$ in (\ref{M_pred}), (\ref{M_q}) and (\ref{M_sigma}). 
Probabilistic upper bounds for the noise and consequences of their 
combination with Theorems \ref{thm-small-lam} and \ref{thm-scaled} 
are discussed in Subsection \ref{prob-bounds}. 
We use the rest of this subsection to discuss $\{M^*_{pred},M^*_q,M^*_\sigma\}$. 

Existing analyses of the Lasso and Dantzig selector can be used find upper bounds for 
$\{M^*_{pred},M^*_q,M^*_\sigma\}$ via 
the sparse eigenvalues 
\cite{CandesT05,CandesT07,ZhangH08,Zhang09-l1,CaiWX10,Zhang10,YeZ10}. 
In the simpler case $A_1=m_1=0$, shaper bounds can be obtained using 
the compatibility factor \cite{vandeGeer07,vandeGeerB09}, 
the restricted eigenvalue \cite{BickelRT09,Koltchinskii09}, or the cone invertibility factors 
\cite{YeZ10,ZhangZ11}. 
Detailed discussions can be found in \cite{vandeGeerB09,YeZ10,ZhangZ11} among others. 
The main difference here is the possibility of excluding some highly correlated vectors from $B$ 
in the case of $A_1>0$. 
The following lemma provide some simple bounds used in our analysis of 
the scaled Lasso estimation of the precision matrix.

\begin{lemma}\label{lm-M^*-bd} 
Let $\{M^*_{pred},M^*_q,M^*_\sigma\}$ be as in (\ref{M_pred}), (\ref{M_q}) and (\ref{M_sigma}) 
with the vector class $\scrU(\bSigmabar,S,B;A,A_1,m,m_1)$ in (\ref{U}).   
Suppose that for a nonnegative-definite matrix $\bSigma$, 
$\max_j\|\bSigmabar_{j,*} - \bSigma_{j,*}\|_\infty\le \lam^*$ and 
$c_*\|\bu_{S\cup B}\|_2^2\le \bu^T\bSigma\bu$ for $\bu\in \scrU(\bSigmabar,S,B;A,A_1,m,m_1)$. 
Suppose further that $\lam^*\{(s_*\vee m)/c_*\}(2A+A_1)^2 \le (A-A_1-1)_+^2/2$. 
Then, 
\bel{lm-M^*-1}
M^*_{pred} + M^*_1\Big(1-\frac{A_1+1}{A}\Big)
\le \max\Big\{\frac{4\vee(4m/s_*)}{c_*(2+A_1/A)^{-2}}, \frac{c_*(1-|S|/s_*)}{A^2}\Big\}
\eel
and
\bel{lm-M^*-2}
M^*_\sigma \le \Big(1+\frac{2A_1}{c_*A}\Big)M^*_{pred}+2(1+A_1)\frac{M^*_1}{A}
 + \frac{A_1m}{As_*} +\frac{2A_1}{A^3}\Big(1-\frac{|S|}{s_*}\Big). 
\eel
Moreover, if in addition $B=\{1,\ldots,p\}$ then 
\bel{lm-M^*-3}
M^*_2  \le (2/c_*)M^*_{pred} + 2(1-|S|/s_*)/(A^2). 
\eel
\end{lemma}

The main condition of Lemma \ref{lm-M^*-bd}, 
\bel{c_*}
c_* \le \inf\left\{\frac{\bu^T\bSigma\bu}{\|\bu_{S\cup B}\|_2^2}: \bu\in \scrU(\bSigmabar,S,B;A,A_1,m,m_1)\right\}, 
\eel
can be viewed as a restricted eigenvalue condition \cite{BickelRT09} on a population version of 
the Gram matrix. 
However, one may also pick the sample version $\bSigma=\bSigmabar$ with $\lam^*=0$. 
Let $\{A,A_1\}$ be fixed constants satisfying $A_1<A-1$. 
Lemma \ref{lm-M^*-bd} asserts that the factors $\{M^*_{pred},M^*_1,M^*_\sigma\}$ 
can be all treated as constants when $1/c_*$ and $m/s_*$ are bounded and 
$\lam^*(s_*\vee m)/c_*$ is smaller than a certain constant. 
Moreover, $M^*_2$ can be also treated as a constant when (\ref{c_*}) holds for $B=\{1,\ldots,p\}$. 

\subsection{Probabilistic error bounds.}\label{prob-bounds}

Theorems \ref{thm-small-lam} and \ref{thm-scaled} provides analytical error bounds
based on the size of the excess noise over a given threshold.
Here we provide probabilistic upper bounds for the excess noise and describe
their implications in combination with Theorems \ref{thm-small-lam} and \ref{thm-scaled}.
We use the following notation:
\bel{L_n}
L_n(t) = n^{-1/2}\Phi^{-1}(1-t),
\eel
where $\Phi^{-1}(t)$ is the standard normal quantile function. 

\begin{proposition}\label{prop-soft-thresh}
Let $\zeta_{(2,m)}(\bv,\lam_*)$ be as in (\ref{zeta})  
and $\kappa_+(m)=\kappa_+(m;\bSigmabar)$ as in (\ref{sparse-eigen}) 
with $\bSigmabar=\bX^T\bX/n$. 
Suppose $\bep\sim N(0,\sigma^2\bI_n)$ and $\|\bx_j\|_2^2=n$. Let $k>0$. \\
(i)  Let $\bz = \bX^T\bep/n$ and $\lam_* = \sigma L_n(k/p)$.
Then, $P\{\zeta_{(2,p)}(\bz,\lam_*)>0\}\le 2k$, and
\bel{prop-soft-thresh-0} \nonumber
& E\zeta_{(2,p)}^2(\bz,\lam_*) \le 4k\lam_*^2/\{L_1^4(k/p)+2L_1^2(k/p)\},
\\ \label{prop-soft-thresh-1} & P\Big\{\zeta_{(2,m)}(\bz,\lam_*) > E\zeta_{(2,p)}(\bz,\lam_*)
+ \sigma L_n(\eps)\sqrt{\kappa_+(m)}\Big\} \le \eps.
\eel
(ii) Let $\sigma^*=\|\bep\|_2/\sqrt{n}$, $\bz^*=\bz/\sigma^*$, 
$\lam_{*,0} = L_{n-3/2}(k/p)$ and $\eps_{n} = e^{1/(4n-6)^2}-1$. Then,
$P\{\zeta_{(2,p)}(\bz^*,\lam_{*,0})>0\}\le (1+\eps_n)k$,
$E\zeta_{(2,p)}^2(\bz^*,\lam_{*,0}) \le (1+\eps_{n})4k\lam_{*,0}^2/\{L_1^4(k/p)+2L_1^2(k/p)\}$, and
\bel{prop-soft-thresh-3}
P\Big\{\zeta_{(2,m)}(\bz^*,\lam_{*,0}) > \mu_{(2,m)}
+L_{n-3/2}(\eps)\sqrt{\kappa_+(m)}\Big\} \le (1+\eps_{n})\eps,
\eel
where $\mu_{(2,m)}$ is the median of $\zeta_{(2,m)}(\bz^*,\lam_{*,0})$. Moreover, 
\bel{prop-soft-thresh-4}
\mu_{(2,m)}\le E\zeta_{(2,p)}(\bz^*,\lam_{*,0}) 
+ (1+\eps_{n})\{\lam_{*,0}/L_1(k/p)\}\sqrt{\kappa_+(m)/(2\pi)}.
\eel
\end{proposition}

We describe consequences of combining Proposition \ref{prop-soft-thresh} with
Theorems \ref{thm-small-lam} and \ref{thm-scaled} in three theorems,
respectively using the probability of no excess noise over the threshold,
the Markov inequality with the second moment, and the concentration bound on the excess noise.

\begin{theorem}\label{cor-regular-lam} Let $0<\eps<p$. Suppose $\bep\sim N(0,\sigma^2\bI_n)$. \\
(i) Let the notation be as in Theorem \ref{thm-small-lam} 
and (\ref{L_n}) with $A_1=0$ and $\lam_* = \sigma L_n(\eps/p^2)$. If (\ref{s_*}) holds, then 
(\ref{th-small-lam-pred-est}) holds with at least probability $1-2\eps/p$. \\
(ii) Let the notation be as in Theorem \ref{thm-scaled} and (\ref{L_n}) with $A_1=0$ 
and $\lam_{*,0} = L_{n-3/2}(\eps/p^2)$.
If (\ref{s_*0-bd}) holds with $P\{(1-\eps_0)^2\le \chi^2_n/n\le (1+\eps_0)^2\}\le \eps/p$, 
then (\ref{thm-scaled-1}) and (\ref{thm-scaled-2}) hold 
with at least probability $1- 3\eps/p$.
\end{theorem}

For a single application of the Lasso or scaled Lasso, 
$\eps/p=o(1)$ guarantees $\|\bz\|_\infty\le \lam_*$ in Theorem \ref{cor-regular-lam} (i) 
and $\|\bz^*\|_\infty\le \lam_{*,0}$ in Theorem \ref{cor-regular-lam} (ii) 
with high probability. The threshold levels are 
$\lam_*/\sigma \approx \lam_{*,0} \approx \lam_{univ} = \sqrt{(2/n)\log p}$, 
as typically considered in the literature. 
In numerical experiments, this often produces nearly optimal results although the 
threshold level may still be somewhat higher than optimal for the prediction and estimation of $\bbeta$. 
However, if we use the union bound to guarantee the simultaneous validity of the oracle inequalities 
in $p$ applications of the scaled Lasso in the estimation of individual columns of a precision matrix, 
Theorem \ref{cor-regular-lam} requires $\eps=o(1)$, or equivalently a significantly higher 
threshold level $\lam_{*,0} \approx \sqrt{(4/n)\log p}$. 
This higher $\lam_{*,0}$, which does not change the theoretical results by much, may
produce clearly suboptimal results in numerical experiments.  

\begin{theorem}\label{cor-small-lam-1}
Let $k>0$. Suppose $\bep\sim N(0,\sigma^2\bI_n)$. \\
(i) Let the notation be as in Theorem \ref{thm-small-lam} and Proposition \ref{prop-soft-thresh}, 
$\lam_* = \sigma L_n(k/p)$, and $A-1>A_1\ge \sqrt{4k/(\eps\, m(L^4_1(k/p)+2L^2_1(k/p)))}$. 
If (\ref{s_*}) holds, then (\ref{th-small-lam-pred-est}) holds with at least 
probability $1-\eps-2|B^c|k/p$. \\
(ii) Let the notation be as in Theorem~\ref{thm-scaled} and Proposition \ref{prop-soft-thresh}, 
{$\lam_{*,0} = L_{n-3/2}(k/p)$, $\eps_{n} = e^{1/(4n-6)^2}-1$,} and $A-1>A_1\ge \sqrt{(1+\eps_n)4k/(\eps\, m(L^4_1(k/p)+2L^2_1(k/p)))}$. 
If (\ref{s_*0-bd}) holds with $P\{(1-\eps_0)^2\le \chi^2_n/n\le (1+\eps_0)^2\}\le \eps$, 
then (\ref{thm-scaled-1}) and (\ref{thm-scaled-2}) hold 
with at least probability $1-2\eps-2|B^c|k/p$.
\end{theorem}

Theorem \ref{cor-small-lam-1} uses the upper bounds for 
$E\zeta_{(2,p)}^2(\bz,\lam_*)$ and $E\zeta_{(2,p)}^2(\bz^*,\lam_{*,0})$ to verify (\ref{noise-cond}). 
Since $L_n(k/p)\approx \sqrt{(2/n)\log(p/k)}$, it allows smaller threshold
levels $\lam_*$ and $\lam_{*,0}$ as long as $k/(\eps\, m(L^4_1(k/p)+2L^2_1(k/p)))$ is small. 
However, it does not allow $\eps\le 1/p$ for using the union bound in $p$ applications of the Lasso 
in precision matrix estimation. 

\begin{theorem}\label{cor-small-lam-2}
Let $k>0$. Suppose $\bep\sim N(0,\sigma^2\bI_n)$. \\
(i) Let the notation be as in Theorem \ref{thm-small-lam} and Proposition \ref{prop-soft-thresh}, 
$\lam_* = \sigma L_n(k/p)$, and
\bes
A-1> A_1\ge \Big(\frac{4k/m}{L^4_1(k/p)+2L^2_1(k/p)}\Big)^{1/2}
+ \frac{L_1(\eps/p)}{L_1(k/p)}\Big(\frac{\kappa_+(m)}{m}\Big)^{1/2}.
\ees
If (\ref{s_*}) holds, then (\ref{th-small-lam-pred-est}) holds with at least probability 
$1-\eps/p-2|B^c|k/p$. \\
(ii) Let the notation be as in Theorem~\ref{thm-scaled} and Proposition \ref{prop-soft-thresh}, 
$\lam_{*,0} = L_{n-3/2}(k/p)$, {$\eps_{n} = e^{1/(4n-6)^2}-1$,} and
\bes
A-1> A_1\ge \Big(\frac{(1+\eps_{n})4k/m}{L^4_1(k/p)+2L^2_1(k/p)}\Big)^{1/2}
+ \Big(\frac{L_1(\eps/p)}{L_1(k/p)}+\frac{1+\eps_{n}}{L_1(k/p)\sqrt{2\pi}}\Big)
\Big(\frac{\kappa_+(m)}{m}\Big)^{1/2}. 
\ees
If (\ref{s_*0-bd}) holds with $P\{(1-\eps_0)^2\le \chi^2_n/n\le (1+\eps_0)^2\}\le\eps/p$, 
then (\ref{thm-scaled-1}) and (\ref{thm-scaled-2}) hold 
with at least probability $1-2\eps/p-2|B^c|k/p$.
\end{theorem}

Theorem \ref{cor-small-lam-2} uses concentration inequalities 
(\ref{prop-soft-thresh-0}) and (\ref{prop-soft-thresh-3}) to verify (\ref{noise-cond}). 
Let $B=\{1,\ldots,p\}$ and $L_n(t)$ be as in (\ref{L_n}). 
By guaranteeing the validity of the oracle inequalities with $1-\eps/p$ probability,
with a reasonably small $\eps$, Theorem \ref{cor-small-lam-2} justifies the use of a fixed
smaller threshold level $\lam_{*,0}=L_{n-3/2}(k/p)\approx \sqrt{(2/n)\log(p/k)}$ 
in $p$ applications of the scaled Lasso to estimate columns of a precision matrix. 

Since $L_1(\eps/p)\approx L_1(k/p)$ typically holds,
Theorem \ref{cor-small-lam-2} only requires $(k/m)/(L^4_1(k/p)+2L^2_1(k/p))$ and $\kappa_+(m)/m$
be smaller than a fixed small constant. This condition relies on the upper sparse eigenvalue only
in a mild way since $\kappa_+(m)/m$ is decreasing in $m$ and
$\kappa_+(m)/m \le 1/m + (1-1/m)\max_{j\neq k}|\bx_j^T\bx_k/n|$ \cite{ZhangH08}.

For $k\asymp m$ and $\log(p/k) \asymp \log(p/(s_*\vee 1))$, 
Theorem \ref{cor-small-lam-2} provides prediction and $\ell_q$ error bounds 
of the orders $\sigma^2(s_*\vee 1)\lam_{*,0}^2\approx \sigma^2((s_*\vee 1)/n)2\log(p/(s_*\vee 1))$ and 
$\sigma (s_*\vee 1)^{1/q}\lam_{*,0}$ respectively. For $\log(p/n)\ll \log n$, 
this could be of smaller order than the error bounds with 
$\lam_{*,0}\approx \lam_{univ}=\sqrt{(2/n)\log p}$. 

Theorem \ref{cor-small-lam-2} suggests the use of a penalty level satisfying 
$\lam/\sigma =\lam_0 = A L_n(k/p)\approx A\sqrt{(2/n)\log(p/k)}$ with
$1< A \le \sqrt{2}$ and a real solution of $k=L_1^4(k/p)+2L_1^2(k/p)$. 
This is conservative since the constraint on $A$ in the theorem is valid with a moderate $m=O(s_*+1)$. 
For $p$ applications of the scaled Lasso in the estimation of precision matrix, this also 
provides a more practical penalty level compared with 
$A'L_n(\eps/p^2)\approx A'\sqrt{(4/n)\log(p/\eps^{1/2})}$, $A'>1$ and $\eps\ll 1$, 
based on existing results and Theorem \ref{cor-regular-lam}.
In our simulation study, we use $\lam_0 = \sqrt{2}L_n(k/p)$ with 
$k=L_1^4(k/p)+2L_1^2(k/p)$.

\subsection{A lower performance bound}
It is well understood that in the class of $\bbeta$ satisfying the sparsity condition (\ref{s_*}), 
$s_*\sigma^2L_n^2(s_*/p)$ and $s_*^{1/q}\sigma L_n(s_*/p)$ are respectively lower bounds 
for the rates of minimax prediction and $\ell_q$ estimation error \cite{YeZ10,RaskuttiWY11}. 
This can be achieved by the Lasso with $\lam_*=\sigma L_n(m/p)$, $m\asymp s_*$, or 
scaled Lasso with $\lam_{*,0}=L_{n-3/2}(m/p)$.  
The following proposition asserts that for each fixed $\bbeta$, the minimax error rate cannot 
be achieved by regularizing the gradient with a threshold level of smaller order. 

\begin{proposition}\label{prop-lower-bd} 
Let $\by=\bX\bbeta+\bep$, $\tbbeta(\lam)$ satisfy 
$\|\bX^T(\by-\bX\tbbeta(\lam))/n\|_\infty\le\lam$, and $\tbh(\lam)=\tbbeta(\lam)-\bbeta$. 
Let $\bSigmabar=\bX^T\bX/n$ and $\kappa_+(m;\cdot)$ be as in (\ref{sparse-eigen}). \\
(i) If $\|\bX^T\bep/n\|_{(2,{k})}\ge {k}^{1/2}\lam_*>0$, then for all $A>1$
\bel{prop-lower-bd-1}
\inf_{\lam\le \lam_*/A} 
\min\Big\{\frac{\|\bX\tbh(\lam)\|^2/n}{\kappa_+^{-1}({{k}};\bSigmabar)},
\frac{\|\tbh(\lam)\|_2^2}{\kappa_+^{-2}({{k}};\bSigmabar)}\Big\}
\ge (1-1/A)^2 k\lam_*^2. 
\eel
(ii) Let $\sigma^*=\|\bep\|_2/\sqrt{n}$ and $N_{k} =\# \{j: |\bx_j^T\bep|/(n\sigma^*) \ge \Ltil_n({k}/p)\}$
with $\Ltil_n(t) = L_n(t)-n^{-1/2}$. 
Suppose $\bX$ has iid $N(0,\bSigma)$ rows, $\diag(\bSigma)=\bI_p$, and 
$2{k} - 4\|\bSigma\|_2 \ge \big(\sqrt{{k}-1}+\sqrt{2\|\bSigma\|_2\log(1/\eps)}\big)^2$. 
Then, $P\{N_{k}\ge {k}\} \ge 1-\eps$ and
\bel{prop-lower-bd-2}
P\Big\{\|\bX^T\bep/n\|_{(q,{k})} 
\ge \sigma^*{k}^{1/q}\sigma\Ltil_n({k}/p)\Big\} \ge 1- \eps. 
\eel
Consequently, there exist numerical constants $c_1$ and $c_2$ such that 
\bes
P\Big\{\inf_{\lam \le c_1\sigma L_n(k/p)} \min\Big(\frac{\|\bX\tbh(\lam)\|^2/n}{\kappa_+^{-1}({{k}};\bSigmabar)},
\frac{\|\tbh(\lam)\|_2^2}{\kappa_+^{-2}({{k}};\bSigmabar)}\Big)
\ge c_2\sigma^2 k L_n^2(k/p)\Big\}\ge 1-\eps - e^{- n/9}. 
\ees
\end{proposition}

It follows from Proposition \ref{prop-lower-bd} (ii) that  
the prediction and $\ell_2$ estimation error is of no smaller order than 
$k\sigma^2L_n^2(k/p)$ for all $\lam \le c_1\sigma L_n(k/p)$. 
This rate is suboptimal when $k \log(p/k) \gg s_*\log (p/s_*)$. 

\section{Estimation after model selection}
We have presented theoretical properties of the scaled Lasso for 
linear regression and precision matrix estimation. 
After model selection, the least squares estimator is often used to remove bias of regularized estimators. 
The usefulness of this technique after the scaled Lasso was demonstrated in \cite{SunZ11}, along with 
its theoretical justification. 
In this section, we extend the theory to smaller threshold level and to the estimation of precision matrix.

In linear regression, the least squares estimator $\bbeta$ and the corresponding estimate of $\sigma$ 
in the model selected by a regularized estimator $\hbbeta$ are given by
\bel{slasso-lse}
\bbetabar=\argmin_{\bb}\Big\{\|y-\bX\bb\|_2^2:\supp(\bb)\subseteq{\widehat S}\Big\},\quad
\sigmabar=\big\|y-\bX\bbetabar\big\|_2\big/\sqrt{n}, 
\eel
where ${\widehat S}=\supp(\hbbeta)$. 
To study the performance of (\ref{slasso-lse}), we define 
sparse eigenvalues relative to a support set $S$ as follows:
\bes
&& \kappa_-^*(m^*,S;\bSigma)
= \min_{J\supseteq S,|{J}\setminus S|\le m^*}
\min_{\|\bu_{J}\|_2=1}\bu_{J}^T\bSigma_{{J},{J}}\bu_{J}, 
\cr && \kappa_+^*(m^*,S;\bSigma) = \min_{{J}\cap S=\emptyset,|{J}|\le m^*} 
\max_{\|\bu_{J}\|_2=1}\bu_{J}^T\bSigma_{{J},{J}}\bu_{J}. 
\ees
It is proved in \cite{SunZ11} that $\{\bbetabar,\sigmabar\}$ satisfies prediction and 
estimation error bounds of the same order as those for the scaled Lasso  
(\ref{scaledlasso}) under some extra conditions on $\kappa_\pm^*(m^*,S;\bSigmabar)$. 
The extra condition on $\kappa_+^*(m^*,S;\bSigmabar)$ is used to derive an 
upper bound for the false positive $|{\widehat S}\setminus S|$, and then the 
extra condition on $\kappa_-^*(m^*,S;\bSigmabar)$ is used to invert $\bX_{S\cup\Shat}$. 
The following theorem extends the result to the smaller threshold level 
$\lam_{*,0}=L_{n-3/2}(k/p)$ in Theorem \ref{cor-small-lam-2} (ii). 
Let 
\bes
M_{lse}^* = \frac{\big[\big\{|S|+\big(\sqrt{m^*}+\sqrt{2m^*\log(ep/m^*)}\big)^2\big\}^{1/2}
+L_1(\eps/p)\big]^2}{s_*\log(p/s_*)}. 
\ees 

\begin{theorem}\label{thm-small-lam-lse}
Let  $(\hbbeta,\hsigma)$ be the scaled lasso estimator in (\ref{scaledlasso}) and
$(\bbetabar,\sigmabar)$ be the least squares estimator (\ref{slasso-lse}) in 
the selected model ${\widehat S}=\supp(\hbbeta)$.
Let the notation
be as in Theorem \ref{cor-small-lam-2} (ii) 
\and $m^* > m$ be an integer satisfying $s_*M^*_{pred}/\{(1-\xi_1)(1-1/A)\}^2 \le m^*/\kappa^*_+(m^*,S)$.
Suppose $\bbeta_{S^c}=0$ and (\ref{noise-cond}) holds with 
$\{\bz_B,\lam_*\}$ replaced by $\{\bz_B^*,\lam_{*,0}\}$. Then, 
\bel{th-mleas-1}
|{\widehat S}\setminus S|\le m^*
\eel
with at least probability $1-2\eps/p-2|B^c|k/p$. Moreover, 
\bel{th-mleas-2}
&& (\sigma^*)^2 - \sigma^2 M_{lse}^*(s_*/n)\log(p/s_*)  
\cr &\le& \sigmabar^2 
\cr &\le& \hsigma^2, 
\cr && \kappa^*_-(m^*-1,S)\|\bhbar\|_2^2 
\\ \nonumber &\le & \|\bX\bhbar\|_2^2/n
\\ \nonumber &\le& (1+\eps_0)^2M_{pred}^*s_*(\sigma\phi_2\lam_0)^2+\sigma^2 M_{lse}^*(s_*/n)\log(p/s_*), 
\eel
with at least probability $1-3\eps/p-2|B^c|k/p$, 
where $\bhbar = \bbetabar-\bbeta$. 
\end{theorem}

Theorem \ref{thm-small-lam-lse} asserts that when $k\vee m^*\asymp s_*$, 
the least squares estimator $\{\bbetabar,\sigmabar\}$ after the scaled Lasso selection 
enjoys estimation and prediction properties 
comparable to that of the scaled Lasso: 
\bes
\lam_0^{-2}\Big\{\big|\sigmabar/\sigma^* -1\big| + \big\|\bbetabar - \bbeta\big\|_2^2 
+ \big\|\bX\bbetabar-\bX\bbeta\big\|_2^2/n\Big\}+\|\bbetabar\|_0 = O_P(1)s_*. 
\ees

Now we apply this method for precision matrix estimation. 
Let $\hbbeta$ be as in (\ref{joint-min}) and define $\bbetabar$ and $\bsigmabar$ as follows: 
\bel{LSE}
&& \bbetabar_{*,j}=\argmin_{\bb}\Big\{\|\bX\bb\|_2^2: b_j=-1,\supp(\bb)\subseteq\supp(\hbbeta_{*,j})\Big\},\ 
\cr && \sigmabar_j =\big\|\bX\bbetabar_{*,j}\big\|_2\big/\sqrt{n}.
\eel
We define $\tbTheta^{\text{LSE}}$ and $\hbTheta^{\text{LSE}}$ as in (\ref{tTheta}) and (\ref{hTheta}) 
with $\bbetabar$ and $\bsigmabar$ in place of $\hbbeta$ and $\hbsigma$.

Under an additional condition on the upper sparse eigenvalue, 
Theorem \ref{thm-small-lam-lse} is parallel to Theorem 8 (ii), 
and the theoretical results in \cite{SunZ11} are parallel to Theorem 6 (ii). 
These results can be used to verify the condition (\ref{oracle-scale}), 
so that Proposition \ref{prop-matrix} also applies to (\ref{LSE}) with the extra upper sparse eigenvalue 
condition on the population correlation matrix $\bR^*$. 
We formally state this result as a corollary. 

\begin{corollary} Under the additional condition $\|\bR^*\|_2=O(1)$ on the population correlation matrix $\bR^*$, 
Theorems \ref{main-1} and \ref{thm-corrmatrix} are applicable to 
the estimator $\hbTheta^{\text{LSE}}$ and the corresponding estimator for $\bOmega^*=(\bR^*)^{-1}$ 
with possibly different numerical constants. 
\end{corollary}

\section{Numerical study}
In this section, we present some numerical comparison between the proposed and existing methods. 
In addition to the proposed estimator (\ref{tTheta}) and (\ref{hTheta})  based on the scaled Lasso (\ref{joint-min}) 
and the least squares estimation after the scale Lasso (\ref{LSE}), 
the graphical Lasso and CLIME are considered. 
The following three models are considered. 
Models 1 and 2 have been considered in \cite{CaiLL11}, while Model 2 in \cite{RothmanBLZ08}.
\begin{itemize}
\item Model 1: $\Theta_{ij}=0.6^{|i-j|}$.
\item Model 2: Let $\bTheta= \bB+\delta\bI$, where each off-diagonal entry in $\bB$ is generated independently 
and equals to 0.5 with probability 0.1 or 0 with probability 0.9. 
The constant $\delta$ is chosen such that the condition number of $\bTheta^*$ is $p$. 
Finally, we rescale the matrix $\bTheta^*$ to the unit in diagonal.
\item Model 3: The diagonal of the target matrix has unequal values. $\bTheta=\bD^{1/2}\bOmega\bD^{1/2}$, 
where $\Omega_{ij}=0.6^{|i-j|}$ and $\bD$ is a diagonal matrix with diagonal elements $d_{ii}=(4i+p-5)/\{5(p-1)\}$.
\end{itemize}
Among the three models, Model 2 is the densest.  
For $p=1000$, the capped $\ell_1$ sparsity $s_*$ is 8.84, 24.63, and 8.80 for three models respectively.

\begin{landscape}
\begin{table}[htbp]\tiny\centering
\caption{Estimation errors under various matrix norms of scaled Lasso, GLasso and CLIME for three models.}\vspace{0.1in}
\begin{tabular}{c rrrr r rrrr r rrrr}\hline\hline\\
\multicolumn{15}{c}{Model 1}\\\\\hline\\
&\multicolumn{4}{c}{Spectrum norm} &&\multicolumn{4}{c}{Matrix $\ell_1$ norm} &&\multicolumn{4}{c}{Frobenius norm} \\ \\\cline{2-5}\cline{7-10}\cline{12-15}\\
$p$ & SLasso & SLasso/LSE & GLasso & CLIME && SLasso & SLasso/LSE & GLasso & CLIME && SLasso & SLasso/LSE & GLasso & CLIME \\\\\hline\\
30  &2.46(0.07) &1.77(0.15)  & 2.49(0.14) & 2.29(0.21) && 2.95(0.10) &2.73(0.21)  & 3.09(0.11) & 2.92(0.17) && 4.20(0.11) &3.37(0.14)& 4.24(0.26) & 3.80(0.36)\\
60  &2.68(0.05) &2.04(0.11)  & 2.94(0.05) & 2.68(0.10) && 3.12(0.08) &3.17(0.25)  & 3.55(0.07) & 3.27(0.09) && 6.41(0.09) &5.35(0.15)& 7.15(0.15) & 6.32(0.28)\\
90  &2.75(0.04) &2.09(0.08)  & 3.07(0.03) & 2.87(0.09) && 3.21(0.07) &3.49(0.31)  & 3.72(0.06) & 3.42(0.07) && 8.09(0.10) &6.87(0.15)& 9.25(0.12) & 8.42(0.31)\\
150 &2.84(0.03) &2.18(0.06)  & 3.19(0.02) & 3.05(0.04) && 3.29(0.07) &3.81(0.31)  & 3.88(0.06) & 3.55(0.06) && 10.79(0.11) &9.32(0.14) & 12.55(0.09) & 11.68(0.20)\\
300 &2.93(0.02) &2.25(0.05)  & 3.29(0.01) & NA         && 3.39(0.05) &4.36(0.38)  & 4.06(0.05) & NA         && 15.83(0.09) &13.89(0.16)& 18.44(0.09) & NA\\
1000&3.08(0.02) &2.38(0.08)  & 3.39(0.00) & NA         && 3.51(0.03) &5.13(0.37)  & 4.44(0.07) & NA         && 30.55(0.09) &26.68(0.19)& 35.11(0.06) & NA\\
\\\hline\hline\\\\
\multicolumn{15}{c}{Model 2}\\\\\hline\\
&\multicolumn{4}{c}{Spectrum norm} &&\multicolumn{4}{c}{Matrix $\ell_1$ norm} &&\multicolumn{4}{c}{Frobenius norm} \\ \\\cline{2-5}\cline{7-10}\cline{12-15}\\
$p$ & SLasso & SLasso/LSE & GLasso & CLIME && SLasso & SLasso/LSE & GLasso & CLIME && SLasso & SLasso/LSE & GLasso & CLIME \\\\\hline\\
30  &0.72(0.08) &1.08(0.17) & 0.82(0.07) & 0.81(0.09) && 1.27(0.15) &1.86(0.33) & 1.49(0.15) & 1.45(0.18) && 1.86(0.10) &2.40(0.24)& 1.84(0.09) & 1.87(0.11)\\
60  &1.06(0.05) &1.41(0.19) & 1.15(0.06) & 1.19(0.08) && 1.93(0.15) &2.68(0.35) & 2.21(0.12) & 2.20(0.23) && 3.27(0.08) &4.19(0.26)& 3.18(0.13) & 3.42(0.09)\\
90  &1.48(0.04) &1.73(0.23) & 1.54(0.05) & 1.61(0.04) && 2.58(0.15) &3.81(0.46) & 2.89(0.16) & 2.90(0.17) && 4.42(0.07) &6.18(0.29)& 4.40(0.11) & 4.65(0.08)\\
150 &1.96(0.03) &2.04(0.28) & 2.02(0.05) & 2.06(0.03) && 3.25(0.17) &5.21(0.63) & 3.60(0.15) & 3.65(0.19) && 5.95(0.06) &9.17(0.33)& 6.19(0.16) & 6.33(0.08)\\
300 &2.88(0.02) &2.34(0.19) & 2.89(0.02) & NA         && 4.45(0.13) &6.75(0.52) & 4.92(0.17) & NA         && 9.26(0.05) &13.99(0.37)& 9.79(0.05) & NA\\
1000&5.46(0.01) &4.88(0.03) & 5.52(0.01) & NA         && 7.09(0.09) &10.26(0.53)& 7.98(0.15) & NA         && 18.85(0.06) &26.08(0.30)& 20.81(0.02) & NA
\\\\\hline\hline\\\\
\multicolumn{15}{c}{Model 3}\\\\\hline\\
&\multicolumn{4}{c}{Spectrum norm} &&\multicolumn{4}{c}{Matrix $\ell_1$ norm} &&\multicolumn{4}{c}{Frobenius norm} \\ \\\cline{2-5}\cline{7-10}\cline{12-15}\\
$p$ & SLasso & SLasso/LSE &GLasso & CLIME && SLasso & SLasso/LSE & GLasso & CLIME && SLasso & SLasso/LSE & GLasso & CLIME \\\\\hline\\
30  & 1.84(0.09) &1.28(0.15) & 2.08(0.10) & 1.63(0.19) && 2.30(0.12) &2.00(0.19) & 2.59(0.10) & 2.17(0.20) && 2.69(0.09) &2.15(0.13)& 2.91(0.16) & 2.37(0.25)\\
60  & 2.18(0.07) &1.58(0.13) & 2.63(0.04) & 2.10(0.10) && 2.63(0.10) &2.46(0.18) & 3.10(0.05) & 2.65(0.14) && 4.10(0.08) &3.41(0.10)& 4.84(0.08) & 3.98(0.13)\\
90  & 2.34(0.06) &1.71(0.11) & 2.84(0.03) & 2.38(0.18) && 2.77(0.10) &2.73(0.20) & 3.30(0.06) & 2.91(0.12) && 5.19(0.08) &4.40(0.10)& 6.25(0.08) & 5.37(0.37)\\
150 & 2.51(0.05) &1.84(0.09) & 3.06(0.02) & 2.76(0.05) && 2.93(0.09) &3.04(0.28) & 3.45(0.04) & 3.18(0.09) && 6.93(0.08) &5.96(0.10)& 8.43(0.07) & 7.75(0.08)\\
300 & 2.70(0.05) &1.99(0.08) & 3.26(0.01) & NA         && 3.10(0.07) &3.39(0.27) & 3.58(0.03) & NA         && 10.18(0.08) &8.89(0.10)& 12.41(0.04) & NA\\
1000& 2.94(0.03) &2.16(0.07) & 3.47(0.01) & NA         && 3.32(0.06) &4.07(0.32) & 3.73(0.03) & NA         && 19.63(0.07) &17.04(0.15)& 23.55(0.02) & NA
\\\\\hline\hline\\
\end{tabular}
\end{table}
\end{landscape}

In each model, we generate a training sample of size 100 from a multivariate normal distribution 
with mean zero and covariance matrix $\bSigma=\bTheta^{-1}$ and an independent sample 
of size 100 from the same distribution for validating the tuning parameter $\lam$ for the graphical 
Lasso and CLIME. 
The GLasso and CLIME estimators are computed based on training data with various $\lam$'s 
and we choose $\lam$ by minimizing likelihood loss 
$\{\text{trace}(\bSigmabar\hbTheta)-\log\det(\hbTheta)\}$ on the validation sample.
The scaled Lasso estimators are computed based on the training sample alone with penalty level 
$\lam_0=AL_n(k/p)$, where $A=\sqrt{2}$ and $k$ is the solution of 
$k=L_1^4(k/p)+2L_1^2(k/p)$. 
The symmetrization step in \cite{CaiLL11} is applied.
We consider six different dimensions $p=30, 60, 90, 150, 300, 1000$ and replicate 100 times in each setting.
The CLIME estimators for $p=300$ and $p=1000$ are not computed due to computational costs.

Table 1 presents the mean and standard deviation of estimation errors based on 100 replications. 
The estimation error is measured by three matrix norms: the spectrum norm, 
the matrix $\ell_1$ norm 
and the Frobenius norm. 
The scaled Lasso estimator, labeled as SLasso, outperforms the graphical Lasso (GLasso) 
in all cases except for the smaller $p\in \{30,60,90\}$ in the Frobenius loss in the denser Model 2. 
It also outperforms the CLIME in most cases, except for smaller $p$ in sparser models 
($p=30$ in Model 1 and $p\in \{30,60\}$ in Model 3).  
The least squares estimator after the scaled Lasso selection outperforms all estimators 
by large margin in the spectrum and Frobenius losses in Models 1 and 3, but in general 
underperforms in the $\ell_1$ operator norm and in Model 2. 
It seems that post processing by the least squares method is a somewhat aggressive procedure 
for bias correction. 
It performs well in sparse models, where variable selection is easier, but may not perform very 
well in denser models. 

Both the scaled Lasso and the CLIME are resulting from sparse linear regression solutions. 
A main advantage of the scaled Lasso over the CLIME is adaptive choice 
of the penalty level for the estimation of each column of the precision matrix. 
The CLIME uses cross-validation to choose a common penalty level for all $p$ columns. 
When $p$ is large, it is computationally difficult. In fact, this prevented us from completing the 
simulation experiment for the CLIME for the larger $p\in\{300,1000\}$.

\section{{Discussion}}

Since the scaled Lasso choose penalty levels adaptively in the estimation 
of each column of the precision matrix, it is expected to outperform methods 
using a fixed penalty level for all columns in the presence of heterogeneity of the diagonal of the 
precision matrix. Let $\tbTheta(\lam)$ be an estimator with columns 
\bel{CLIME}
{\tbTheta}_{*j}(\lam) = \argmin_{\bv\in\R^p}\Big\{\bigl\|\bv\bigr\|_{1}: 
\bigl\| {\bSigmabar}\bv - \bfe_{j}\bigr\|_{\infty}\leq \lambda\Big\},\ j=1,\ldots, p. 
\eel
The CLIME is a symmetrization of this estimator $\tbTheta(\lam)$ with fixed 
penalty level for all columns. 
In the following example, the scaled Lasso estimator has a faster convergence rate 
than (\ref{CLIME}). The example also demonstrates the possibility of achieving 
the rate $d\lam_0$ in Theorem \ref{main-1} with unbounded $\|\bTheta^*\|_2\ge d^2$, 
when Theorem \ref{main} is not applicable. 

\begin{example}\label{example-1} 
Let $p > n^2 + 3+m$ with ${(m, m^4(\log p)/n)} \to (\infty,{0)}$ and  $4m^2\le \log p$. 
Let $L_n(t)\approx \sqrt{(2/n)\log(1/t)}$ be as in (\ref{L_n}). 
Let $\{J_1,J_2,J_3\}$ be a partition of $\{1,\ldots,p\}$ with $J_1=\{1,2\}$ 
and $J_2=\{3,\ldots,3+m\}$. Let $\rho_1 = \sqrt{1-1/m^2}$, 
$\bv = (v_1,\ldots,v_{m})^T \in\R^{m}$ with $v_j^2=1/m$, $\rho_2=c_0 m^{3/2} L_n(m/p)=o(1)$, and 
\bes
\bSigma^*  &= &\begin{pmatrix}
\bSigma^*_{J_1,J_1} &0&0 \cr 
0 & \bSigma^*_{J_2,J_2} & 0 \cr 0 & 0 & \bI_{p-m-3} 
\end{pmatrix},\ 
\bSigma^*_{J_1,J_1} = \begin{pmatrix} 1 & \rho_1 \cr \rho_1  & 1 \end{pmatrix},\ 
\bSigma^*_{J_2,J_2} = \begin{pmatrix} 1 & \rho_2 \bv^T \cr \rho_2 \bv  & \bI_{m} \end{pmatrix}. 
\ees
The eigenvalues of $\bSigma^*_{J_1,J_1}$ are $1\pm\rho_1$, those of 
$\bSigma^*_{J_2,J_2}$ are $1\pm\rho_2,1,\ldots,1$, and
\bes
(\bSigma^*_{J_1,J_1})^{-1} = m^2\begin{pmatrix} 1 & - \rho_1
\cr  - \rho_1 & 1 \end{pmatrix},\ 
(\bSigma^*_{J_2,J_2})^{-1} = \frac{1}{1-\rho_2^2}\begin{pmatrix} 1 & - \rho_2 \bv^T 
\cr  - \rho_2 \bv & (1-\rho_2^2)\bI_{m}+\rho_2^2\bv\bv^T \end{pmatrix}. 
\ees
We note that $\diag(\bSigma^*)=\bI_p$, $d=m+1$ is the maximum degree,  
$\|\bTheta^*\|_2=1/(1-\rho_1)\approx 2d^2$, and $\|\bTheta^*\|_1 \approx 2d^2$. 
The following statements are proved in the Appendix. \\
(i) Let $\hbTheta$ be the scaled Lasso estimator of $\bTheta^*=(\bSigma^*)^{-1}$ with penalty level 
$\lam_0=A\sqrt{(4/n)\log p}$, $A>1$, as in Theorem \ref{main-1}.  
Then, there exists a constant $M^*_1$ such that 
\bes
P\Big\{ \|\hbTheta-\bTheta^*\|_2\le \|\hbTheta-\bTheta^*\|_1\le M^*_1 m L_n(m/p)\Big\} \to 1. 
\ees
(ii) If $\rho_2=c_0 m^{3/2} L_n(m/p)$ with a sufficiently small constant $c_0>0$, then 
\bes
P\Big\{ \inf_{\lam>0}\|\tbTheta(\lam)-\bTheta^*\|_2
\ge c_0 m^{3/2} L_n(m/p)/\sqrt{1+1/m}\Big\}\to 1. 
\ees
Thus, the order of the $\ell_1$ and spectrum norms of the error of (\ref{CLIME}) for the best data dependent 
penalty level $\lam$ is larger than that of the scaled Lasso 
by a factor $\sqrt{m}$.
\end{example}

\section{Proofs}
In this section, we provide all proofs. 
We first prove the results in Section 4 since they are used to prove the results in Section 3. 

\paragraph{{Proof of Proposition \ref{prop-sparse-norm}}.} 
Lemma 20 in \cite{YeZ10} gives $\|\bv\|_q^q\le \|\bv\|_{(q,{m})}^q+(a_q/{m})^{q-1}\|\bv\|_1^q$. 
The rest of part (i) follows directly from definition. 
Lemma 20 in \cite{YeZ10} also gives 
$\|\bv\|_{(q,{m})}^* \le \|\bv\|_{(q',{m}/a_q)} + {m}^{-1/q}\|\bv\|_1$. 
The rest of part (ii) is dual to the corresponding parts of part (i). Since 
$\|\bSigmabar\bv\|_{(2,{m})} =\max_{\|\bu\|_0=m,\|\bu\|_2=1}\bu^T\bSigmabar\bv$ 
and $\|\bu^T\bSigmabar^q\|_2\le \kappa_+^q({m};\bSigmabar)$ for $q\in \{1/2,1\}$, 
part (iii) follows. 
$\hfill\square$ 

\paragraph{{Proof of Theorem \ref{thm-small-lam}}.} 
By the Karush-Kuhn-Tucker conditions, 
\bel{KKT-1}
(\bX^T\bX/n)\bh = \bz - \lam\bg,\ \sgn(\hbeta_j)g_j\in \{0,1\},\
\|\bg\|_\infty\le 1.
\eel
Since $\zeta_{(2,m)}(\bz_B,\lam_*)$ is the $\|\cdot\|_{(2,m)}$ norm of $(|\bz_B|-\lam_*)_+$, 
(\ref{noise-cond}) implies 
\bel{h-z-1}
|\bh^T\bz| &\le& \lam_*\|\bh\|_1 + \sum_{j\in B} |h_j|(|z_j|-\lam_*)_+ 
\cr &\le& \lam_*\|\bh\|_1 + A_1\lam_*m^{1/2}\|\bh_B\|_{(2,m)}^*. 
\eel
Since $ - h_j\sgn(\hbeta_j) \le |\beta_j|- |\hbeta_j|\le \min(|h_j|,- |h_j|+2|\beta_j|)\ \forall\ j\in S^c$, 
(\ref{s_*}) and (\ref{KKT-1}) yield
\bes
- \lam \bh^T \bg
&\le& \lam \|\bh_S\|_1 - \lam \|\bh_{S^c}\|_1 + 2\lam \|\bbeta_{S^c}\|_1
\cr &\le& \lam \|\bh_S\|_1 - \lam \|\bh_{S^c}\|_1 + 2\lam\lam_*(s_*-|S|). 
\ees
By applying the above bounds to the inner product of $\bh$ and (\ref{KKT-1}), we find 
\bes
\|\bX\bh\|_2^2/n
&\le& A_1\lam_*m^{1/2}\|\bh_B\|_{(2,m)}^* + (\lam_*- \lam) \|\bh_{S^c}\|_1 
\cr && + (\lam+\lam_*) \|\bh_S\|_1 + 2\lam\lam_*(s_*-|S|). 
\ees
Let $\bu = (A/\lam)\bh$. It follows that when $A\lam_* \le \lam$, 
\bes
\frac{\|\bX\bu\|_2^2}{n}
\le A_1m^{1/2}\|\bu_B\|_{(2,m)}^*
- (A-1)\|\bu_{S^c}\|_1 + (A+1)\|\bu_S\|_1 + 2A(s_*-|S|). 
\ees
Since $\bX^T\bX/n=\bSigmabar$, $\bu\in \scrU(\bSigmabar,S,B;A,A_1,m,m_1)$ 
with $m_1=s_*-|S|$. 
Since $\bh  = \lam\bu/A$ and $s_*=m_1+|S|$, 
the conclusion follows from (\ref{M_pred}) and (\ref{M_q}). $\hfill\square$

\paragraph{{Proof of Theorem \ref{thm-scaled}}.} 
It follows from the scale equivariance of (\ref{scaledlasso}) that
\bel{scaledlasso-1}
\{\hbbeta/\sigma^*,\hsigma/\sigma^*\} =
\{\hbb,\hphi\}
= \argmin_{\bb,\phi}\Big\{\|\by^*-\bX\bb\|_2^2/(2n\phi)+\lam_0\|\bb\|_1+\phi/2\Big\},
\eel
where $\by^* = \by/\sigma^*= \bX\bb^* + \bep^*$ with
$\bb^* = \bbeta/\sigma^*$ and $\bep^*=\bep/\sigma^*$. 
Our objective is to bound $\|\bX(\hbb-\bb^*)\|_2^2/n$ and $\|\hbb-\bb^*\|_q$ 
from the above and $\hsigma/\sigma^*$ from both sides.
To this end, we apply Theorem \ref{thm-small-lam} to the Lasso estimator
\bes
\hbb(\lam) = \argmin_{\bb}\Big\{\|\by^*-\bX\bb\|_2^2/(2n)+\lam\|\bb\|_1\Big\}.
\ees

Let $\bz^*=\bz/\sigma^*$ and $\bh^*(\lam) = \hbb(\lam) - \bb^*$. 
Since $\|\by^*-\bX\bb^*\|_2^2/n=\|\bep^*\|_2^2/n=1$, 
\bes
1-\|\by^*-\bX\hbb(\lam)\|_2^2/n
&=& \bh^*(\lam)^T\bX^T(\by^*-\bX\hbb(\lam))/n + \bh^*(\lam)^T\bz^*
\cr &=& 2\bh^*(\lam)^T\bz^* - \|\bX\bh^*(\lam)\|_2^2/n. 
\ees
Consider $\lam \ge A\lam_{*,0}$.  
Since (\ref{noise-cond}) holds with 
$\{\bz,\lam_*\}$ replaced by $\{\bz^*,\lam_{*,0}\}$, we find as in the proof of 
Theorem \ref{thm-small-lam} that  
\bes
\bu(\lam) &=& \bh^*(\lam)A/\lam\in  \scrU(\bSigmabar,S,B;A,A_1,m,m_1). 
\ees
In particular, (\ref{h-z-1}) gives 
\bes 
|\bh^*(\lam)^T\bz^*| &\le& \lam_{*,0}\|\bh^*(\lam)\|_1+A_1\lam_{*,0}m^{1/2}\|\bh^*_B(\lam)\|_{(2,m)}^*
\cr &\le & (\lam^2/A^2)\Big\{\|\bu(\lam)\|_1+A_1m^{1/2}\|\bu_B(\lam)\|_{(2,m)}^*\Big\}. 
\ees
Thus, the definition of $M^*_\sigma$ in (\ref{M_sigma}) gives 
\bes
|2\bh^*(\lam)^T\bz^* - \|\bX\bh^*(\lam)\|_2^2/n| < M^*_\sigma s_*\lam^2. 
\ees
We summarize the calculation in this paragraph with the following statement: 
\bel{basic-scaled}
\lam \ge A\lam_{*,0}\ \Rightarrow\ 
\big|1-\|\by-\bX\hbb(\lam)\|_2^2/n\big| < M^*_\sigma s_*\lam^2.
\eel

As in \cite{SunZ11}, 
the convexity of the joint loss function in (\ref{scaledlasso-1}) implies 
\bes
(\phi  - \hphi)(\phi^2  - \|\by-\bX\hbb(\phi\lam_0)\|_2^2/n) \ge 0,
\ees 
so that $\hphi$ can be bounded by testing the sign of $\phi^2  - \|\by-\bX\hbb(\phi\lam_0)\|_2^2/n$. 
For $(\phi,\lam)=(\phi_1,\phi_1\lam_0)$, we have 
\bes
\lam^2 
= \frac{\lam_0^2}{1+\lam_0^2M^*_\sigma s_*}
\ge \frac{A^2\lam_{*,0}^2}{1-\eta_*+A^2\lam_{*,0}^2 M^*_\sigma s_*} = A^2\lam_{*,0}^2, 
\ees
which implies $\|\by-\bX\hbb(\phi_1\lam_0)\|_2^2/n > 1-\phi_1^2\eta_0 = \phi_1^2$ by (\ref{basic-scaled}) 
and the definition of $\phi_1$. This yields $\hphi > \phi_1$.  Similarly, $\hphi < \phi_2$. 
The error bounds for the prediction and the estimation $\hbbeta$
follow from Theorem \ref{thm-small-lam} due to 
$A\lam_{*,0}\le \phi_1\lam_0 < \hphi\lam_0 < \phi_2\lam_0$.
$\hfill\square$

\paragraph{{Proof of Lemma \ref{lm-M^*-bd}}.} 
By Proposition \ref{prop-sparse-norm}, 
$m^{1/2}\|\bu_B\|_{(2,m)}^*\le \|\bu_B\|_1+m^{1/2}\|\bu_B\|_{(2,4m)}$, 
so that for $\bu\in \scrU(\bSigmabar,S,B;A,A_1,m,m_1)$, 
\bes
\bu^T\bSigmabar\bu+(A-A_1-1)\|\bu\|_1
\le 2A\|\bu_S\|_1+A_1m^{1/2}\|\bu_B\|_2  + 2A m_1. 
\ees
Let $\xi=A/(A-A_1-1)$ and $\xi_1=A_1/(A-A_1-1)$. It follows that 
\bel{pf-lm-1-1}
&& (\xi/A)\bu^T\bSigmabar\bu + \|\bu\|_1 
\cr &\le& 2\xi\|\bu_S\|_1+\xi_1m^{1/2}\|\bu_B\|_{2}+ 2\xi m_1
\cr &\le&(2\xi |S|+\xi_1 m+2\xi m_1)^{1/2}\{(2\xi+\xi_1)\|\bu_{S\cup B}\|_2^2+ 2\xi m_1\}^{1/2}
\cr &\le&\{(2\xi s_*+\xi_1 m)/c_*\}^{1/2}\{(2\xi+\xi_1)\bu^T\bSigma\bu+2\xi c_* m_1\}^{1/2}
\cr &\le&\{(s_*\vee m)/c_*\}^{1/2}(2\xi+\xi_1)(\bu^T\bSigma\bu+c_* m_1)^{1/2}
\eel
due to $s_*=|S|+m_1$ and $c_*\|\bu_{S\cup B}\|_2^2\le \bu^T\bSigma\bu$. 
In terms of $\{\xi,\xi_1\}$, the condition of the Lemma can be stated as 
$\lam^*\{(s_*\vee m)/c_*\}(2\xi+\xi_1)^2 \le 1/2$. Thus, 
\bel{pf-lm-1-2}
\bu^T\bSigma\bu - \bu^T\bSigmabar\bu \le \lam^*\|\bu\|_1^2
\le \bu^T\bSigma\bu/2 + c_* m_1/2.
\eel 
Inserting this inequality back into (\ref{pf-lm-1-1}), we find that 
\bes
(\xi/A)\bu^T\bSigmabar\bu + \|\bu\|_1
\le \{(s_*\vee m)/c_*\}^{1/2}(2\xi+\xi_1)(2\bu^T\bSigmabar\bu+2c_* m_1)^{1/2}.
\ees
If $(\xi/A)\bu^T\bSigmabar\bu + \|\bu\|_1\ge (\xi/A)(2\bu^T\bSigmabar\bu+2c_*m_1)/4$, 
we have 
\bes
(\xi/A)\bu^T\bSigmabar\bu + \|\bu\|_1 
\le \{(s_*\vee m)/c_*\}(2\xi+\xi_1)^2(4A/\xi). 
\ees
Otherwise, we have 
$(\xi/A)\bu^T\bSigmabar\bu + 2\|\bu\|_1\le (\xi/A)c_*m_1$. Consequently, 
\bes
(\xi/A)\bu^T\bSigmabar\bu + \|\bu\|_1
\le \max\Big\{\{(s_*\vee m)/c_*\}(2\xi+\xi_1)^2(4A/\xi),(\xi/A)c_*(s_*-|S|)\Big\}. 
\ees
This and the definition of $\{M^*_{pred},M^*_1\}$ yield (\ref{lm-M^*-1}) via 
\bes
\xi M^*_{pred} + M^*_1
\le \max\Big\{\big(1\vee(m/s_*))(2\xi+\xi_1)^2\frac{4}{\xi c_*}, \xi c_*(1-|S|/s_*)/A^2\Big\}. 
\ees
Moreover, (\ref{pf-lm-1-2}) gives 
$c_*\|\bu_{S\cup B}\|_2^2 \le \bu^T\bSigma\bu \le 2\bu^T\bSigmabar\bu +2 c_*m_1$, so that $M^*_\sigma$ 
can be bounded via 
\bes
&& \bu^T\bSigmabar\bu/(s_*A^2)+2\big(\|\bu\|_1+A_1m^{1/2}\|\bu_B\|_{(2,m)}^*\big)/(s_*A^2)
\cr &\le&  M^*_{pred}+2(1+A_1)\|\bu\|_1/(s_*A^2)+ (A_1/A)\Big\{m/s_*+\|\bu_B\|_2^2/(s_*A^2)\Big\}
\cr &\le&  M^*_{pred}+2(1+A_1)M^*_1/A+(A_1/A)\Big(m/s_* + (2/c_*)M^*_{pred}+2(1-|S|/s_*)/A^2\Big). 
\ees
This gives (\ref{lm-M^*-2}). 
If in addition $B=\{1,\ldots,p\}$, then 
\bes
M^*_2 = \sup_{\bu\in\scrU}\|\bu\|_2^2/(s_*A^2)
\le (2/c_*)M^*_{pred} + 2(1-|S|/s_*)/A^2. 
\ees
This completes the proof.  $\hfill\square$

\bigskip
The tail probability bound for $\zeta^*(\bz^*,\lam_*,m)/\sigma^*$ in part (ii) of 
Proposition \ref{prop-soft-thresh} uses the following version of the L\'evy concentration inequality in the sphere.

\begin{lemma}\label{lm-sphere}
Let $\teps_m = \sqrt{2/(m-1/2)}\Gamma(m/2+1/2)/\Gamma(m/2)-1$,
${\bU}=({U}_1,\ldots,{U}_{m+1})^T$ be a uniform random vector in $S^m=\{\bu\in\R^{m+1}:\|\bu\|_2=1\}$,
$f(\bu)$ a unit Lipschitz function in $S^m$,
and $m_f$ the median of $f({\bU})$. Then,
\bel{lm-sphere-1}
P\{ {U}_1 > x \} \le (1+\teps_m) P\big\{ N(0,1/(m-1/2)) > \sqrt{-\log(1-x^2)}\big\},
\eel
$1 < 1+\teps_m  < \exp\big(1/(4m-2)^{2}\big)$, and
\bel{lm-sphere-2}
P\{ f({\bU}) > m_f + x\} &\le& P\Big\{{U}_1 > x\sqrt{1-(x/2)^2}\Big\}
\cr &\le& (1+\teps_m)P\{N(0,1/(m-1/2))>x\}.
\eel
\end{lemma}

{\bf Proof.} Since ${U}_1^2$ follows the beta$(1/2,m/2)$ distribution,
\bes
P\big\{ {U}_1 > x \big\} = \frac{\Gamma(m/2+1/2)/2}{\Gamma(m/2)\Gamma(1/2)}
\int_{x^2}^1 t^{-1/2}(1-t)^{m/2-1}dt.
\ees
Let $y = \sqrt{-(m-1/2)\log(1-t)}$.
We observe that $-t^{-1}\log(1-t)\le (1-t)^{-1/2}$ by inspecting the infinite series expansions
of the two functions. This gives
\bes
\frac{e^{-y^2/2}dy}{t^{-1/2}(1-t)^{m/2-1}dt} = \frac{t^{1/2}e^{-y^2/2}(m-1/2)^{1/2}}{2(-\log(1-t))^{1/2}(1-t)^{m/2}}
\ge 2^{-1}\sqrt{m-1/2}.
\ees
Since $y=\sqrt{-(m-1/2)\log(1-x^2)}$ when $t=x^2$, it follows that
\bes
P\big\{ {U}_1 > x \big\} \le (1+\teps_m)
\int_{\sqrt{-(m-1/2)\log(1-x^2)}}^\infty (2\pi)^{-1/2}e^{-y^2/2}dt.
\ees
Let ${\Atil} = \{\bu\in S^m: f(\bu) \le m_f\}$, $H=\{\bu: u_1\le 0\}$, and
$A_x =\cup_{\bv\in A}\{\bu\in S^m: \|\bu-\bv\|_2\le x\}$ for all $A\subset S^m$.
Since $\bu\in {\Atil}_x$ implies $f(\bu) \le m_f + x$ and $P\{\bU\in {\Atil}\}\ge P\{\bU\in H\}$,
the L\'evy concentration inequality gives
\bes
P\{ f({\bU}) > m_f + x\} \le P\{{\bU}\not\in H_x\} = P\Big\{{U}_1 > x\sqrt{1-(x/2)^2}\Big\}.
\ees
The second inequality of (\ref{lm-sphere-2}) then follows from
$(d/dx)\{- \log\{1-(x^2-x^4/4)\} - x^2\}\ge 0$ for $x^2\le 2$ and $\|{U}_1\|_\infty\le 1$.

It remains to bound $1+\teps_m$. Let $x=m+1/2$. Since
\bes
\frac{1+\teps_m}{1+\teps_{m+2}} = \frac{(m/2)\sqrt{m+3/2}}{(m/2+1/2)\sqrt{m-1/2}}
= \frac{(x-1/2)\sqrt{x+1}}{(x+1/2)\sqrt{x-1}},
\ees
the infinite series expansion of its logarithm is bounded by
\bes
\log\Big(\frac{1+\teps_m}{1+\teps_{m+2}}\Big)
= \frac{1}{2}\log\Big(\frac{1+1/x}{1-1/x}\Big)+\log\Big(\frac{1-1/(2x)}{1+1/(2x)}\Big)
\le \frac{x^{-3}}{4} + \frac{x^{-5}}{5}+\cdots
\ees
Since $\{(x-1)^{-2}-(x+1)^{-2}\}/2
= 2 x^{-3} + 4x^{-5} + \cdots $ by Newton's binomial formula,
\bes
\log\Big(\frac{1+\teps_m}{1+\teps_{m+2}}\Big)
\le \{(x-1)^{-2}-(x+1)^{-2}\}/16.
\ees
This gives $\log(1+\teps_m) \le 1/\{16(x-1)^2\}$.
$\hfill\square$

\paragraph{{Proof of Proposition \ref{prop-soft-thresh}}.} (i)
Let $L=L_1(k/p)$. Since $P\{ N(0,\sigma^2/n) > \lam_*\} = k/p$,
$\lam_*=\sigma L/\sqrt{n}$.
Since $z_j = \bx_j^T\bep/n\sim N(0,\sigma^2/n)$, $P\{\zeta_{(2,p)}(\bz,\lam_*)>0\}\le 2k$ and 
\bes
E\zeta_{(2,p)}^2(\bz,\lam_*) &=& p(\sigma^2/n)E(|N(0,1)|-L)_+^2
\cr &=& 2p(\sigma^2/n)\int_{L}^\infty (x-L)^2\varphi(x)dx.
\ees
Let $J_k(t)=\int_0^\infty x^ke^{-x-x^2/(2t^2)}dx$. By definition
\bes
\frac{t^2\int_t^\infty (x-t)^2\varphi(x)dx}{\Phi(-t)}
=\frac{t^2\int_0^\infty x^2e^{-t x -x^2/2}dx}{ \int_0^\infty e^{-t x -x^2/2}dx}
=\frac{\int_0^\infty u^2e^{- u -u^2/(2t^2)}du}{ \int_0^\infty e^{- u -u^2/(2t^2)}du}
=\frac{J_2(t)}{J_0(t)}.
\ees
Since $J_{k+1}+J_{k+2}/t^2 = - \int_0^\infty x^{k+1} de^{-x-x^2/(2t^2)} = (k+1)J_k(t)$, we find
\bes
\frac{J_2(t)}{J_0(t)}
=\frac{J_2(t)}{\{J_2(t)+J_3(t)/t^2\}/2+J_2(t)/t^2} \le \frac{1}{1/2+1/t^2}.
\ees
Thus, $E\zeta_{(2,p)}^2(\bz,\lam_*) = 2p(\sigma^2/n)(k/p)L^{-2}J_2(L)/J_0(L)
\le 2k\lam_*^2 L^{-4}/(1/2+1/L^2)$.

Since $\bz_j = \bx_j^T\bep/n$, $\big(\sum_{j\in B}(|z_j| - \lam_*)_+^2\big)^{1/2}$ is a
function of $\bep$ with the Lipschitz norm $\|\bX_B/n\|_{2}$.
Thus, $\zeta_{(2,m)}(\bz,\lam_*)$ is a function of $\bep$ with the  Lipschitz norm
$\max_{|B|=m}\|\bX_B/n\|_2=\sqrt{\kappa_+(m)/n}$.
In addition, since $\zeta_{(2,m)}(\bz,\lam_*)$ is an increasing convex function of $(|z_j|-\lam_*)_+$ 
and $(|z_j|-\lam_*)_+$ are convex in $\bep$, $\zeta_{(2,m)}(\bz,\lam_*)$ is  
a convex function of $\bep$. The mean of $\zeta_{(2,m)}(\bz,\lam_*)$ is no smaller than its median.
This gives (\ref{prop-soft-thresh-1}) by the Gaussian concentration inequality \cite{Borell75}.

(ii) The scaled version of the proof uses Lemma \ref{lm-sphere} with $m=n-1$ there.
Let ${\bU}=\bep/\|\bep\|_2$,
$z^*_j = \bx_j^T\bep/(n\sigma^*) = (\bx_j/\sqrt{n})^T{\bU}$ and $\bz^* = \bX^T\bep/(n\sigma^*)$.
Since $z^*_j\sim {U}_1$, (\ref{lm-sphere-1}) yields the bound 
$P\{\zeta_{(2,p)}(\bz^*,\lam_{*,0})>0\}\le (1+\eps_n)2k$ and 
\bes
E\zeta_{(2,p)}^2(\bz^*,\lam_{*,0})
= p E(|{U}_1|-\lam_*)_+^2 \le (1+\eps_{n})p E(|N(0,1)|-L)_+^2/(n-3/2). 
\ees
The bound for $E\zeta_{(2,p)}^2(\bz^*,\lam_{*,0})$ is then derived as in (i). 
Lemma \ref{lm-sphere} also gives
\bes
&& P\big\{\pm(\zeta_{(2,m)}(\bz^*,\lam_{*,0}) - \mu_{(2,m)}) 
>  x \sqrt{\kappa_+(m)}\big\} 
\cr &\le& (1+\eps_{n})P\{|N(0,1/(n-3/2))|>x\}
\ees
and
\bes
|E\zeta_{(2,m)}(\bz^*,\lam_{*,0})-\mu_{(2,m)}| 
&\le& (1+\eps_{n})\sqrt{\kappa_+(m)/(n-3/2)}E\big(N(0,1)\big)_+
\cr &=& (1+\eps_{n})\sqrt{\kappa_+(m)/\{2\pi(n-3/2)\}}.
\ees
The above two inequalities yield (\ref{prop-soft-thresh-3}) and (\ref{prop-soft-thresh-4}). $\hfill\square$

\paragraph{Proof of Theorems \ref{cor-regular-lam}, \ref{cor-small-lam-1} and \ref{cor-small-lam-2}.} 
The conclusions follow from Theorems \ref{thm-small-lam} and~\ref{thm-scaled} once 
(\ref{noise-cond}) is proved to hold with the given probability. 
In Theorem \ref{cor-regular-lam}, the tail probability bounds for $\zeta_{(0,p)}$ in 
Proposition \ref{prop-soft-thresh} yield (\ref{noise-cond}) with $A_1=0$. 
In Theorem~\ref{cor-small-lam-1}, the moment bounds for $\zeta_{(0,p)}$ in 
Proposition \ref{prop-soft-thresh} controls the excess noise in (\ref{noise-cond}). 
In Theorem \ref{cor-small-lam-2} (i), we need 
$A_1\lam_*m^{1/2} \ge E\zeta_{(0,p)}(\bz,\lam_*)+\sigma L_n(\eps/p)\sqrt{\kappa_+(m)}$
by (\ref{prop-soft-thresh-0}), so that the given lower bound of $A_1$ suffices due to 
$L_n(\eps)/L_n(k/p)=L_1(\eps/p)/L_1(k/p)$. 
The proof of Theorem \ref{cor-small-lam-2} (ii) is nearly identical, with (\ref{prop-soft-thresh-3}) 
and (\ref{prop-soft-thresh-4}) in place of (\ref{prop-soft-thresh-3}). 
We omit the details. 
$\hfill\square$

\paragraph{{Proof of Proposition \ref{prop-lower-bd}}.} 
(i) By the $\ell_\infty$ constraint, 
\bes
\|\bX^T(\bep - \bX\tbh(\lam))/n\|_{(2,{k})} = \|\bX^T(\by - \bX\tbbeta(\lam))/n\|_{(2,{k})} \le \lam\sqrt{{k}}.
\ees 
Thus, when $\lam\sqrt{{k}} \le \|\bX^T\bep/n\|_{(2,{k})}/A$, 
\bes
\|\bX^T\bep/n\|_{(2,{k})}(1-1/A) \le \|\bX^T\bep/n\|_{(2,{k})} - \lam\sqrt{{k}} 
\le \|\bX^T\bX\tbh(\lam)/n\|_{(2,{k})}. 
\ees
Thus, Proposition \ref{prop-sparse-norm} (iii) gives (\ref{prop-lower-bd-1}).

(ii) Let $f(x) = (x - \Ltil_1({k}/p))_+\wedge 1$ and $\bz^* = \bX^T\bep/\|\bep\|_2$. 
Since $\bz^*\sim N(0,\bSigma)$ and $\|f(\bz^*)\|_2$ has unit Lipschitz norm, 
the Gaussian concentration theorem gives 
\bes
P\Big\{ Ef(\bz^*) - f(\bz^*) \ge \sqrt{2\|\bSigma\|_2\log(1/\eps)}\Big\} \le \eps. 
\ees
This implies $\Var(f(\bz^*))\le 4\|\bSigma\|_2$. 
Since $Ef^2(\bz^*) \ge p P\{|N(0,1)|\ge L_1({k}/p)\} = 2{k}$, 
\bes
Ef(\bz^*) - \sqrt{2\|\bSigma\|_2\log(1/\eps)}
\ge \sqrt{2{k} - 4\|\bSigma\|_2} - \sqrt{2\|\bSigma\|_2\log(1/\eps)} 
\ge \sqrt{{k}-1}. 
\ees
This gives $P\{N_{k}\ge {k}\} \ge P\{f(\bz^*) > \sqrt{{k}-1}\} \ge1-\eps$ due to $N_{k}\ge f^2(\bz^*)$. 
Thus, (\ref{prop-lower-bd-2}) follows from 
 $\|\bX^T\bep/n\|_{(q,{k})} \ge \sigma^*{k}^{1/q}\Ltil_n({k}/p)$ when $N_{k}\ge {k}$. The final conclusion 
 follows from part (i) and large deviation for $(\sigma^*/\sigma)^2\sim \chi^2_n/n$. 
$\hfill\square$ 

\begin{lemma}\label{lm-4} Let $\chi^2_{m,j}$ be $\chi^2$ distributed variables with $m$ degrees of freedom. Then, 
\bes
E\max_{1\le j\le t}\chi^2_{m,j} 
\le \big(\sqrt{m} + \sqrt{2\log  t}\big)^2, \quad t\ge 1. 
\ees
\end{lemma}

{\bf Proof.} Let $f(t) = 2\log t - \int_0^\infty \min\big(1, t P\{N(0,1)> x\}\big) dx^2$. 
We first proof $f(t)\ge 0$ for $t\ge 2$.  Let $L_1(x)=-\Phi^{-1}(x)$. 
We have $f(2) \ge 2\log 2 - 1>0$ and 
\bes
f'(t) = 2/t - 2\int_{L_1(1/t)}^\infty P\big\{N(0,1)>x\big\}xdx 
\ge 2/t-2\int_{L_1(1/t)}^\infty \varphi(x) dx = 0. 
\ees
The conclusion follows from 
$P\{\chi_{m,j} > \sqrt{m}+x\}\le P\{N(0,1)>x\}$ for $x>0$. $\hfill\square$ 

\paragraph{Proof of Theorem \ref{thm-small-lam-lse}.} Let $\bh=\hbbeta-\bbeta$ and $\hlam=\hsigma\lam_0$.
Consider ${J}\subseteq {\widehat S}\setminus S$ with $m\le |{J}|\le m^*$. 
For any $j\in {\widehat S}$, it follows from the KKT conditions that 
$|\bx_j^T\bX\bh/n|=|\bx_j^T(\by-\bX\hbbeta-\bep)|\ge \hlam-|z_j|$.
By the definition of $\kappa^*_+(m^*,S)$ and (\ref{th-small-lam-pred-est}),
\bel{pf-thmlse-1}
\sum_{j\in {J}} (\hlam-|z_j|)_+^2 
&\le & \sum_{j\in {J}}|\bx_j^T\bX\bh/n|^2
\cr &=&(\bX_{J}^T\bX\bh/n)^T(\bX_{J}^T\bX\bh/n)
\cr &\le& \kappa^*_+(m^*,S)\|\bX\bh\|_2^2/n\cr
&\le& \kappa^*_+(m^*,S)M^*_{pred}s_*\hlam^2. 
\eel
Since $\zeta_{(2,k)}(\bz^*_B,\lam_{*,0})/k^{1/2}\downarrow k$ by Proposition \ref{prop-sparse-norm} (i), 
the $\{z^*,\lam_{0,*}\}$ version of (\ref{noise-cond}) gives 
$\zeta_{(2,|J|)}(\bz^*,\lam_{*,0})/|J|^{1/2}\le \zeta_{(2,m)}(\bz^*_B,\lam_{*,0})/m^{1/2}\le\xi_1(A-1)\lam_{*,0}$. 
Thus, with $z^*_j=z_j/\sigma^*$, 
\bes
\sum_{j\in {J}} (\hlam-|z_j|)_+^2 
&\ge& \sum_{j\in {J}} \Big\{\hlam-\sigma^*\lam_{*,0}-\sigma^*(|z_j^*|-\lam_{*,0})_+\Big\}_+^2
\cr &\ge& \Big\{|{J}|^{1/2}(\hlam-\sigma^*\lam_{*,0})-\sigma^*\zeta_{(2,|{J}|)}(\bz^*_B,\lam_{*,0})\Big\}_+^2
\cr &\ge& |J| \Big\{\hlam-\sigma^*\lam_{*,0} -\sigma^*\xi_1(A-1)\lam_{*,0}\Big\}_+^2
\ees
Since $\lam_{*,0}^2/(\lam_0\phi_1)^2 = (\lam_{*,0}/\lam_0)^2(1+\eta_0)
\le (1-\eta_*)/A^2+\eta_*/A^2 =1/A^2$, 
we have $\lam_{*,0}\sigma^* < \lam_{*,0}\hsigma/\phi_1 \le \hlam/A$. 
The above inequalities and (\ref{pf-thmlse-1}) yield
\bes
|J|\le \frac{\kappa^*_+(m^*,S)M^*_{pred}s_*\hlam^2}
{\big\{\hlam-\sigma^*\lam_{*,0} -\sigma^*\xi_1(A-1)\lam_{*,0}\big\}_+^2}
< \frac{\kappa^*_+(m^*,S)M^*_{pred}s_*}
{\big\{1-1/A -\xi_1(1-1/A)\big\}_+^2}\le m^*. 
\ees
Since ${\widehat S}\setminus S$ does not have a subset of size $m^*$, 
we have $|{\widehat S}\setminus S|<m^*$ as stated in (\ref{th-mleas-1}). 
Let $\bP_B$ be the projection to the linear span of $\{\bx_j, j\in B\}$. We have 
\bel{pf-thm-lse-1}
&& \hsigma^2 \ge \|\bP_{\Shat}^\perp\by\|_2^2/n = \sigmabar^2 
\ge \|\bP_{S\cup\Shat}^\perp\by\|_2^2/n 
= (\sigma^*)^2 - \|\bP_{S\cup\Shat}\bep\|_2^2/n,
\cr && \|\bX\bhbar\|_2^2/n = \|\bP_{\Shat}\by-\bX\bbeta\|_2^2/n
= \|\bP_{\Shat}\bep\|_2^2/n+\|\bP_{\Shat}^\perp\bX\bbeta\|_2^2/n. 
\eel

Let $N = {p\choose m^*}$. We have $\log N\le m^*\log(ep/m^*)$ by Stirling. By Lemma \ref{lm-4}, 
\bes
E\|\bP_{S\cup\Shat}\bep\|_2^2/\sigma^2
\le E\max_{|B|=m^*}\|\bP_{S\cup B}\bep\|_2^2/\sigma^2
\le |S|+(\sqrt{m^*}+\sqrt{2\log N})^2. 
\ees
Since $\max_{|B|=m}\|\bP_{S\cup B}\bep\|_2$ is a unit Lipschitz function, 
\bes
\|\bP_{S\cup\Shat}\bep\|_2/\sigma
&\le& \big\{|S|+\big(\sqrt{m^*}+\sqrt{2m^*\log(ep/m^*)}\big)^2\big\}^{1/2} + L_1(\eps/p)
\cr &\le& \sqrt{M_{lse}^*s_*\log(p/s_*)}
\ees
with probability $\eps/p$. 
In addition, Theorem \ref{thm-scaled} gives $\|\bP_{\Shat}^\perp\bX\bbeta\|_2^2\le \|\bX\hbbeta-\bX\bbeta\|_2^2
\le (1+\eps_0)^2M_{pred}^*s_*(\sigma\phi_2\lam_0)^2$. 
Inserting these bounds into (\ref{pf-thm-lse-1}) yields (\ref{th-mleas-2}).  
$\hfill\square$

{\lemma\label{chisq} 
Suppose that the rows of $\bX\in \mathbb{R}^{n\times p}$ are iid $N(0,\bSigma)$ random vectors. \\
(i) Let $Y=\trace(\bA\bX'\bX/n)$ and 
$\sigma^2=\trace\{(\bA+\bA')\bSigma(\bA+\bA')\bSigma\}/2$ with a deterministic matrix $\bA$. 
Then, $EY=\mu=\trace(\bA\bSigma)$, $\Var(Y)=\sigma^2/n$ and 
\bes
E\exp\Big\{t(Y-\mu)\Big\}\le \exp\Big\{-\frac{t\sigma}{\sqrt{2}}-\frac{n}{2}\log(1-\sqrt{2}t\sigma/n)\Big\}. 
\ees 
Consequently, for $0<x\le 1$, 
\bes
P\Big\{(Y-\mu)/\sigma>x\Big\}\le \exp\Big\{-\frac{n}{2}\Big(\sqrt{2}x-\log(1+\sqrt{2}x)\Big)\Big\}\le e^{-nx^2/4}. 
\ees
(ii) Let $\bR^*$ and $\overline{\bR}$ be the population and sample correlation matrices of $\bX$. Then, 
\bes
P\Big\{|\overline{R}_{jk}-R^*_{jk}| > x\sqrt{1-(R^*_{jk})^2}\Big\}\le 2 P\big\{|t_n|>n^{1/2}x\Big\}
\ees
where $t_n$ has the t-distribution with $n$ degrees of freedom. 
In particular, for $n \ge 4$, 
\bes
P\Big\{|\overline{R}_{jk}-R^*_{jk}| > \sqrt{2}x\Big\}\le 2e^{1/(4n-2)^2}  P\big\{ |N(0,1/n)| > x\big\},\ 0\le x\le 1. 
\ees
}

{\bf Proof of Lemma \ref{chisq}.} 
(i) This part can be proved by computing the moment generating function 
with $t\sigma/n =x/(1+\sqrt{2}x)$. We omit details. For $0<x<1$, 
\bes
f(x) = \frac{\sqrt{2}x - \log(1+\sqrt{2}x)}{x^2} 
= \int_0^{\sqrt{2}x} \frac{u du}{x^2(1+u)}
= \int_0^{\sqrt{2}} \frac{u du}{1+xu}\ge f(1)>1/2. 
\ees
(ii) Conditionally on $\Sigmabar_{kk}$, 
$\Sigmabar_{jk}/\Sigmabar_{kk} \sim N(\Sigma_{jk}/\Sigma_{kk}, (1-(R^*_{jk})^2)\Sigma_{jj}/(n\Sigmabar_{kk}))$. Thus, 
\bes
z_{jk} 
&=& \Big(\frac{n\Sigmabar_{kk}}{(1-(R^*_{jk})^2)\Sigma_{jj}}\Big)^{1/2}
\Big(\frac{\Sigmabar_{jk}}{\Sigmabar_{kk}} - \frac{\Sigma_{jk}}{\Sigma_{kk}}\Big)
\cr &=& \Big(\frac{n}{1-(R^*_{jk})^2}\Big)^{1/2}
\Big(\overline{R}_{jk}\frac{\Sigmabar_{jj}^{1/2}}{\Sigma_{jj}^{1/2}} 
- R_{jk}\frac{\Sigmabar_{kk}^{1/2}}{\Sigma_{kk}^{1/2}}\Big)
\ees
is a $N(0,1)$ variable independent of $\Sigmabar_{kk}$. Consequently, 
\bes
\frac{n^{1/2}|\overline{R}_{jk}-R_{jk}|}{(1-(R^*_{jk})^2)^{1/2}}
=\frac{|z_{jk}+z_{kj}|}{\Sigmabar_{jj}^{1/2}/\Sigma_{jj}^{1/2}+\Sigmabar_{kk}^{1/2}/\Sigma_{kk}^{1/2}}
\le |t_{jk}|\vee|t_{kj}|
\ees
with $t_{jk} = z_{jk}\Sigma_{kk}^{1/2}/\Sigmabar_{kk}^{1/2}\sim t_n$. 
Let $U_1$ be a uniformly distributed variable in the unit sphere of $\R^{n+1}$. 
Since $t_n^2/n \sim U_1^2/(1-U_1^2)$, Lemma \ref{lm-sphere} provides 
\bes
P\big\{ t_n^2/n > e^{x^2}-1\big\} 
= P\big\{ U_1^2 > 1- e^{-x^2}\big\} \le 2e^{1/(4n-2)^2}P\Big\{N(0,1/(n-1/2)) > x\Big\}. 
\ees
The conclusion follows from 
$e^{x^2n/(n-1/2)}-1\le 2x^2$ for $0<x\le 1$. 
$\hfill\square$

\medskip
{\bf Proof of Proposition \ref{prop-matrix}.} 
Since $\Thetatil_{jj} = 1/\hsigma_j^2$ and $\max\{C_0\lam_0, C_1s_{*,j}\lam_0^2\}\le 1/4$, 
(\ref{oracle-scale}) and the condition on $\sigma^*_j$ implies
\bes
|\Thetatil_{jj}/\Theta^*_{jj}|\le (5/4)^3\le 2,\ 
|\Thetatil_{jj}/\Theta^*_{jj}-1|\le \{(5/4)^2C_0+5/4+1\}\lam_0. 
\ees 
It follows from (\ref{tTheta}), (\ref{oracle-scale}) and the condition on 
${\overline \bD}=\diag(\Sigmabar_{jj},j\le p)$ that 
\bes
&& \big\|\tbTheta_{*,j}-\bTheta^*_{*,j}\big\|_1 =
\big\| - \hbbeta_{*,j}\Thetatil_{jj} -\bTheta^*_{*,j}\big\|_1 
\cr && \quad\qquad \le \big\| (\hbbeta_{-j,j} -\bbeta_{-j,j})\Thetatil_{jj}\big\|_1 
+\big\|\bTheta^*_{*,j}(\Thetatil_{jj}/\Theta^*_{jj}-1)\big\|_1  
\cr && \quad\qquad \le 
\big\|\hbD_{-j}^{-1/2}\big\|_\infty \big\|\hbD_{-j}^{1/2}(\hbbeta_{-j,j} - \bbeta_{-j,j})\big\|_1\Thetatil_{jj} 
+\big\|\bTheta^*_{*,j}\big\|_1\big|\Thetatil_{jj}/\Theta^*_{jj}-1\big| 
\cr && \quad\qquad \le 
(5/2)\Theta^*_{jj}\big\|\bD_{-j}^{-1/2}\big\|_\infty
(\Theta^*_{jj})^{-1/2}C_2 s_{*,j}\lam_0+\big\|\bTheta_{*,j}\big\|_1\{(3/2)C_0+5/2\}\lam_0
\cr && \quad\qquad \le
C\Big\{(\big\|\bD_{-j}^{-1}\big\|_\infty\Theta^*_{jj})^{1/2}s_{*,j}\lam_0 +\big\|\bTheta_{*,j}\big\|_1\lam_0\Big\}
\ees
with $C = \max(5C_2/2,3C_0/2+5/2)$. This gives (\ref{bound-hbTheta}) due to 
$\|\hbTheta-\bTheta^*\|_1\le 2\|\tbTheta-\bTheta^*\|_1$ by (\ref{hTheta}). 
Similarly, 
\bes
\big\|\tbOmega_{*,j}-\bOmega^*_{*,j}\big\|_1 
&=& 
\big\| - \hbD_{-j}^{1/2}\hbbeta_{*,j}\Thetatil_{jj}\Dhat_{jj}^{1/2} -\bOmega^*_{*,j}\big\|_1 
\cr &\le & \big\|\hbD_{-j}^{1/2}(\hbbeta_{-j,j} - \bbeta_{-j,j})\big\|_1 \Thetatil_{jj}\Dhat_{jj}^{1/2} 
\cr && 
+\big\|\hbD_{-j}^{1/2}\bD_{-j}^{-1/2}\bOmega^*_{*,j}(\Thetatil_{jj}/\Theta^*_{jj})(\Dhat_{jj}/D_j)^{1/2} 
-\bOmega^*_{*,j}\big\|_1   
\cr &\le & C\big\{(\Theta^*_{jj})^{-1/2}s_{*,j}\lam_0 \Theta^*_{jj}D_{jj}^{1/2}
+\big\|\bOmega_{*,j}\big\|_1\lam_0\big\}
\ees
This gives (\ref{bound-hbOmega}) due to $D_{jj}\Theta^*_{jj} = \Omega^*_{jj}$. 
We omit an explicit calculation of $C$. 

Let $\chi_{n,j}^2 = n\bTheta_{jj}^*(\sigma^*_j)^2$. When $\chi_{n,j}^2 \sim \chi_n^2$, we have 
\bes
|\Thetatil_{jj}/\Theta^*_{jj}-1|\le \{(5/4)^2+5/4\}C_1s_{*,j}\lam_0^2+(4/3)|\chi^2_{n,j}/n-1|
\ees
It follows from Lemma \ref{chisq} that 
$P\{ |\chi^2_{n,j}/n-1| > \sqrt{2}x\} \le 2e^{- nx^2/4}$ for $x\le 1$. 
Let $a_j = \|\bTheta_{*,j}\|_1$, $t = \max\{M\max_ja_j/\sqrt{n},\tau_n(\bTheta^*)\}$ 
and $B_0=\{j: a_j\le \sqrt{8}t\}$. By definition $t\le M\tau_n(\bTheta^*)$ and 
$nt^2/a_j^2\ge M^2$. It follows that 
\bes
P\Big\{\max_j |\chi^2_{n,j}/n-1|a_j > 4t \Big\}
&\le& |B_0|e^{-n/4} + \sum_{j\not\in B_0}\exp(-2n t^2/a_j^2)
\cr &\le& pe^{-n/4} + e^{-M^2}\sum_{j\not\in B_0}\exp\Big(- n \tau^2(\bTheta^*)/a_j^2\Big). 
\cr &\le& pe^{-n/4} + e^{-M^2}. 
\ees
Thus, $\max_j |\Thetatil_{jj}/\Theta^*_{jj}-1|a_j = O_P\big(\tau_n(\bTheta^*)+\max_j s_{*,j}a_j\lam_0^2)$. 
$\hfill\square$ 

\medskip
{\bf Proof of Theorem \ref{main-1}.}  
We need to verify conditions (\ref{oracle-scale}) and (\ref{prop-matrix-1}) in order to apply 
Proposition \ref{prop-matrix}. 
Since $\Theta_{jj}^*(\sigma^*_j)^2\sim \Sigmabar_{jj}/\Sigma^*_{jj}\sim \chi^2_n/n$, 
(\ref{prop-matrix-1}) follows from Lemma \ref{chisq}~(i) with $\lam_0\asymp \sqrt{(\log p)/n}$. 
Moreover, the condition $P\{(1-\eps_0)^2\le\chi^2_n/n\le (1+\eps_0)^2\}\le\eps/p$ holds with 
small $\eps_0$ and $\eps$ since $\sqrt{(\log p)/n} = \lam_0/(2A)$ is 
assumed to be sufficiently small. We take $\eps_0=0$ in (\ref{s-starj}) since its value 
does not change the order of $s_{*,j}$. 

If we treat $\Sigmabar_{kk}^{1/2}\beta_k$ as the regression coefficient in (\ref{joint-min}) 
for the standardized design vector $\Sigmabar_{kk}^{-1/2}\bx_k$, $k\neq j$, 
Theorem \ref{cor-regular-lam} (ii) asserts that the conclusions of Theorem \ref{thm-scaled} 
hold with probability $1-3\eps/p$ for each $j$, with $\lam_0=A\sqrt{4(\log p)/n}$, 
$A_1=0$ and $\eps\asymp 1/\sqrt{\log p}$. 
By the union bound, the conclusions of Theorem \ref{thm-scaled} holds simultaneously for all $j$ with probability 
$1-3\eps$. Moreover, (\ref{oracle-scale}) is included in the conclusions of Theorem \ref{thm-scaled} 
when $M^*_\sigma$ and $M^*_1$ are uniformly bounded in the $p$ regression problems 
with large probability. Thus, it suffices to verify the uniform boundedness of these quantities. 

We use Lemma \ref{lm-M^*-bd} to verify the uniform boundedness of $M^*_\sigma$ and $M^*_1$ 
with $A_1=0$, $B_j=S_j$, $m_j=0$ and $\{\bSigmabar,\bSigma^*\}$ replaced by 
$\{{\overline \bR}_{-j,-j},\bR^*_{-j,-j}\}$. 
Note that the Gram matrix for the regression problem in (\ref{joint-min}) is 
${\overline \bR}_{-j,-j}$, which is random and dependent on $j$, so that 
$M^*_\sigma$ and $M^*_{1}$ are random and dependent on $j$ with the random design. 
It follows from Lemma \ref{chisq} (ii) that 
\bes
\max_{k\neq j}\| {\overline \bR}_{k,-j}-\bR^*_{k,-j}\|_\infty
\le \max_{j,k}|{\overline R}_{k,j}-R_{k,j}|\le L_n(5\eps/p^2)
\ees
with probability $1-\eps$. We may take $L_n(5\eps/p^2) = 2\sqrt{(\log p)/n}$ with $\eps\asymp 1/\sqrt{\log p}$. 
This yields the first condition of Lemma \ref{lm-M^*-bd} with $\lam^*=2\sqrt{(\log p)/n}\asymp \lam_0$.  
The second condition $c_*\|\bu_S\|_2^2\le \bu^T\bR^*_{-j,-j}\bu$ 
follows from (\ref{partial-inv}). The third condition translates to 
$\max_{j\le p}\lam_0s_{*,j}\le c_0$, which is imposed in Theorem \ref{main-1}. 
Thus, all conditions of Lemma \ref{lm-M^*-bd} hold simultaneously for all $j$ with large probability. 
The proof is complete since the conclusions of Lemma~\ref{lm-M^*-bd} with $m=m_j=0$ 
guarantee the uniform boundedness of $M^*_\sigma$ and $M^*_{1}$. $\hfill\square$ 

\paragraph{{Proof of Theorem \ref{thm-corrmatrix}.}} 
The proof is parallel to that of Theorem \ref{main-1}. 
Since the smaller $\lam_{*,0}=L_{n-3/2}(k/p)$ is used, 
we need to apply Theorem \ref{cor-small-lam-2} (ii) with  
$A_1>0$, $m=m_j>0$ and typically much larger $B_j$ than $S_j$. 
Since the condition $m_j\le C_0s_{*,j}$ is impose in (\ref{cond-xi1}), 
the conclusions of Lemma \ref{lm-M^*-bd} still guarantee 
the uniform boundedness of $M^*_\sigma$ and $M^*_{1}$. 
The verification of the conditions of Lemma \ref{lm-M^*-bd} is identical to 
the case of larger $\lam_{*,0}$ in Theorem \ref{main-1}. 
The only difference is the need to verify that condition (\ref{cond-xi1}) uniformly 
guarantees the condition on $A_1$ in Theorem \ref{cor-small-lam-2} (ii), 
where $\kappa_+/m$ has the interpretation of 
$\kappa_+(m_j;{\overline \bR}_{-j,-j})/m_j$, which depends on $j$ and 
random ${\overline \bR}$. 
Anyway, it suffices to verify $\kappa_+(m_j;{\overline \bR}_{-j,-j})/m_j\le\psi_j$ 
simultaneously for all $j$ with large probability.  

We verify $\kappa_+(m_j;{\overline \bR}_{-j,-j})/m_j\le\psi_j$ with the same 
argument as in Lemma \ref{lm-M^*-bd}. For any vector $\bu$ with $\|\bu\|_0=m_j$ 
and $\|\bu\|_2=1$, it holds with probability $1-\eps$ that 
\bes
\Big|\bu^T({\overline \bR}_{-j,-j}-\bR_{-j,-j})\bu\Big|
\le \max_{j,k}|{\overline R}_{k,j}-R_{k,j}|\sum_{j,k}|u_ju_k|
\le L_n(5\eps/p^2) m_j. 
\ees
Thus, it follows from the definition of $\kappa_+(m;\bSigma)$ in (\ref{sparse-eigen}) that 
$\kappa_+(m_j;{\overline \bR}_{-j,-j})/m_j
\le \kappa_+(m_j;{\bR}_{-j,-j})/m_j + L_n(5\eps/p^2)=\psi_j$ for all $j$. This completes the proof. 
$\hfill\square$

\paragraph{{Proof of Example \ref{example-1}}.}
(i) Let $s_{*,j}=d_j=\#\{k: \Theta_{jk}^*\neq 0\}\le m+1$. We have $\max_j(1+s_{*,j})\lam_0
\le (m+2)\lam_0\to 0$. 
Let $B_j=\{k\neq j:\Theta^*_{kj}\neq 0\}$. 
Since $B_j=J_1\setminus\{j\}$ for $j\in J_1$, (\ref{partial-inv}) holds with 
\bes
\inf\Big\{\bu^T(\bR^*_{-j,-j})\bu/\|\bu_{B_j}\|_2^2: \bu_{B_j}\neq 0\Big\} \ge 1-\rho_2\to 1. 
\ees
Thus, Theorem \ref{main-1} is directly applicable to this example. 

Next, we calculate the error bound in (\ref{bound-hbTheta}) and (\ref{main-1-3}). 
Since $d_j(\bTheta^*_{jj})^{1/2}=2/(1-\rho_1^2)^{1/2}=2m$ for $j\in J_1$ and 
$d_j(\bTheta^*_{jj})^{1/2}\le (m+1)/(1-\rho_2^2)^{1/2}\le 2m$ for $j\in J_2$, 
\bes
(\big\|\bD_{-j}^{-1}\|_\infty\Theta^*_{jj})^{1/2}s_{*,j}\lam_0
= (\Theta^*_{jj})^{1/2}s_{*,j}\lam_0 \le 2m\lam_0. 
\ees
In addition, $\|\bTheta_{*,j}\|_1 \le 2m^2$ for $j\in J_1$ and 
$\|\bTheta_{*,j}\|_1 \le (1+\rho_2\|\bv\|_1)/(1-\rho_2^2)\le 3/2+o(1)$ 
for $j\in J_2$, so that for $t = \sqrt{(2/n)\log p}$, 
\bes
\sum_j \exp(-n t^2/\|\bTheta_{*,j}\|_1^2)
\le 2\exp\Big(-\frac{2\log p}{4m^2}\Big) + p\exp\Big(-\frac{2\log p}{3/2+o(1)}\Big) \to 0. 
\ees 
It follows that the quantities in (\ref{main-1-3}) are bounded by 
\bes
\max_{j\le p}s_{*,j}\|\bTheta^*_{*,j}\|_1\lam_0^2\le 2(m\lam_0)^2,\ 
\tau_n(\bTheta^*)\le \sqrt{(2/n)\log p}\le \lam_0/(A\sqrt{2}). 
\ees
Since $m\lam_0\to 0$, the error for the scaled Lasso is of the order 
$m\lam_0$ by Theorem \ref{main-1}. 
The conclusion follows since $L_n(m/p)=(1+o(1))\sqrt{(2/n)\log p}$ 
when $4m^2\le\log p$. 

(ii)Let $\tlam = \max_j\|\bSigmabar_{*,j}-\bSigma^*_{*,j}\|_\infty$ and $\tlam^*=\rho_2/\sqrt{m}+\tlam$. 
Since $\diag(\bSigma^*)=\bI_n$, 
\bes
\tlam \lesssim L_n(1/p) \ll \rho_2/\sqrt{m},\quad 
\tlam^* = (1+o(1))\rho_2/\sqrt{m}=(1+o(1))c_0mL_n(m/p). 
\ees 
For $\lam \ge \tlam^*$, $\bfe_3$ is feasible for (\ref{CLIME}) with $j =3 \in J_2$, 
so that $\| {\tbTheta}_{*j}(\lam)\|_1\le 1$. Since $\|\bTheta^*_{J_2,3}\|_1\ge 1+m^{1/2}\rho_2$, 
\bes
m^{1/2}\rho_2 \le \inf_{\lam\ge\tlam^*}\|{\tbTheta}_{J_2,3}(\lam) - {\bTheta}^*_{J_2,3}(\lam)\|_1
\le (m+1)^{1/2}\inf_{\lam\ge\tlam^*}\|{\tbTheta}_{J_2,3}(\lam) - {\bTheta}^*_{J_2,3}(\lam)\|_2. 
\ees
It follows that for $\lam\ge\tlam^*$, ${\tbTheta}(\lam)$ is suboptimal in the sense of 
\bes
\inf_{\lam\ge\tlam^*}\|{\tbTheta}(\lam) - {\bTheta}^*(\lam)\|_2 \ge {\sqrt{m/(1+m)}}\rho_2
=c_0 m^{3/2} L_n(m/p)/\sqrt{1+1/m}. 
\ees

Consider $\lam\le\tlam^*$. Let $\bbeta_{-j,j}=- {\bTheta}^*_{-j,j}/{\bTheta}^*_{jj}$, 
$\tbbeta_{-j,j}(\lam)=- {\tbTheta}_{-j,j}(\lam)/{\tbTheta}_{jj}(\lam)$, $\sigma_j=(\bTheta^*_{j,j})^{-1/2}$, 
and $\tbh_j(\lam) = \tbbeta_{-j,j}(\lam) - \bbeta_{-j,j}$. 
By (\ref{CLIME}),
$\|\bX_{-j}^T(\bx_j - \bX_{-j}\tbbeta(\lam))/n\|_\infty \le \lam/\tbTheta_{jj}(\lam)$. 
Since $m(\log p)/n\to 0$ and $\|\bSigma^*\|_2\le 2$, $P\{\kappa_+(m;\bSigmabar)\le 3\}\to 1$. 
Thus, by Proposition \ref{prop-lower-bd}, there exist positive constants $\{c_1,c_2\}$ such that 
\bes
\min_j P\Big\{\inf_{\lam/\tbTheta_{jj}(\lam)\le c_1\sigma_jL_n(m/p)}
\|\tbh_j(\lam)\|_2 \ge c_2\sigma_j \sqrt{m} L_n(m/p)\Big\}\to 1. 
\ees
For $\tbTheta_{jj}(\lam) \ge \bTheta^*_{jj}/2$, 
\bes
\|\tbh_j(\lam)\|_2
&=& \|\tbTheta_{-j,j}(\tlam)/\tbTheta_{jj}(\lam) - \bTheta^*_{-j,j}/\bTheta^*_{jj}\|_2
\cr &\le&\|\tbTheta_{-j,j}(\tlam)- \bTheta^*_{-j,j}\|_2/\tbTheta_{jj}(\lam) 
+ \|\bbeta_{-j,j}\|_2|\tbTheta_{jj}(\lam)-\bTheta^*_{jj}|/\tbTheta_{jj}(\lam)
\cr &\le& \|\tbTheta_{*,j}(\tlam)- \bTheta^*_{*,j}\|_2(1+\|\bbeta_{-j,j}\|_2)/(\bTheta^*_{jj}/2). 
\ees
For $j=1$, $\bTheta^*_{jj}=m^2$ and $\|\bbeta_{-j,j}\|_2=\rho_1$, 
so that $\|\tbTheta_{*,j}(\tlam)- \bTheta^*_{*,j}\|_2\ge m^2\|\tbh_j(\lam)\|_2/4$
when $\tbTheta_{jj}(\lam) \ge \bTheta^*_{jj}/2$. 
Since $\|\tbTheta_{*,j}(\tlam)- \bTheta^*_{*,j}\|_2\ge m^2/2$ when $\tbTheta_{jj}(\lam) \le \bTheta^*_{jj}/2$, 
\bes
\inf_{\lam\le\tlam^*} \|\tbTheta(\lam)-\bTheta^*\|_2\ge \min\Big(m^2\|\tbh_1(\lam)\|_2/4, m^2/2\Big). 
\ees
Pick $0< c_0 < \min(c_1/2,c_2/4)$. Since $\sigma_1=(\bTheta^*_{11})^{-1/2}=1/m$, 
\bes
&& P\Big\{\inf_{\lam\le\tlam^*} \|\tbTheta(\lam)-\bTheta^*\|_2
\le \min\big(m^2/2,(c_2/4)m^{3/2} L_n(m/p)\big)\Big\}
\cr &\le & P\Big\{\tlam^* > (m^2/2)(c_1/m)L_n(m/p)\Big\}+o(1) = o(1).  
\ees
Since $L_n(m/p)\to 0$ implies 
$\min\big(m^2/2,(c_2/4)m^{3/2} L_n(m/p)\big) \ge c_0m^{3/2}L_n(m/p)$, the conclusion follows. $\hfill\square$

\bibliographystyle{abbrv}
\bibliography{MatrixEst}

\begin{thebibliography}{10}

\bibitem{AbramovichG10}
F.~Abramovich and V.~Grinshtein.
\newblock Map model selection in gaussian regression.
\newblock {\em Electr. J. Statist.}, 4:932--949, 2010.

\bibitem{Antoniadis10}
A.~Antoniadis.
\newblock Comments on: $\ell_1$-penalization for mixture regression models.
\newblock {\em Test}, 19(2):257--258, 2010.

\bibitem{BanerjeeEA08}
O.~Banerjee, L.~El~Ghaoui, and A.~d'Aspremont.
\newblock Model selection through sparse maximum likelihood estimation for
  multivariate gaussian or binary data.
\newblock {\em Journal of Machine Learning Research}, 9:485--516, 2008.

\bibitem{BickelRT09}
P.~Bickel, Y.~Ritov, and A.~Tsybakov.
\newblock Simultaneous analysis of {L}asso and {D}antzig selector.
\newblock {\em Annals of Statistics}, 37(4):1705--1732, 2009.

\bibitem{BirgeM01}
L.~Birge and P.~Massart.
\newblock Gaussian model selection.
\newblock {\em J. Eur. Math. Soc.}, 3:203--268, 2001.

\bibitem{BirgeM07}
L.~Birge and P.~Massart.
\newblock Minimal penalties for gaussian model selection.
\newblock {\em Probability Theory Related Fields}, 138:33--73, 2007.

\bibitem{Borell75}
C.~Borell.
\newblock The brunn-minkowski inequality in gaussian space.
\newblock {\em Invent. Math.}, 30:207--216, 1975.

\bibitem{BuneaTW07-Aggr}
F.~Bunea, A.~Tsybakov, and M.~Wegkamp.
\newblock Aggregation for gaussian regression.
\newblock {\em Ann. Statist.}, 35:1674--1697, 2007.

\bibitem{CaiLL11}
T.~Cai, W.~Liu, and X.~Luo.
\newblock A constrained $\ell_1$ minimization approach to sparse precision
  matrix estimation.
\newblock {\em Journal of the American Statistical Association}, 106:594--607,
  2011.

\bibitem{CaiWX10}
T.~Cai, L.~Wang, and G.~Xu.
\newblock Shifting inequality and recovery of sparse signals.
\newblock {\em IEEE Transactions on Signal Processing}, 58:1300--1308, 2010.

\bibitem{CandesT05}
E.~J. Cand{\`e}s and T.~Tao.
\newblock Decoding by linear programming.
\newblock {\em IEEE Trans. on Information Theory}, 51:4203--4215, 2005.

\bibitem{CandesT07}
E.~J. Cand{\`e}s and T.~Tao.
\newblock The dantzig selector: statistical estimation when $p$ is much larger
  than $n$ (with discussion).
\newblock {\em Annals of Statistics}, 35:2313--2404, 2007.

\bibitem{DonohoJ94}
D.~L. Donoho and I.~Johnstone.
\newblock Minimax risk over $\ell_p$--balls for $\ell_q$--error.
\newblock {\em Probability Theory and Related Fields}, 99:277--303, 1994.

\bibitem{FriedmanHT08}
J.~Friedman, T.~Hastie, and R.~Tibshirani.
\newblock Sparse inverse covariance estimation with the graphical {L}asso.
\newblock {\em Biostatistics}, 9:432--441, 2008.

\bibitem{HuberR09}
P.~J. Huber and E.~M. Ronchetti.
\newblock {\em Robust statistics}, pages 172--175.
\newblock Wiley, second edition, 2009.

\bibitem{Koltchinskii09}
V.~Koltchinskii.
\newblock The dantzig selector and sparsity oracle inequalities.
\newblock {\em Bernoulli}, 15:799--828, 2009.

\bibitem{LamF09}
C.~Lam and J.~Fan.
\newblock Sparsistency and rates of convergence in large covariance matrices
  estimation.
\newblock {\em Annals of Statistics}, 37:4254--4278, 2009.

\bibitem{MeinshausenB06}
N.~Meinshausen and P.~B{\"u}hlmann.
\newblock High-dimensional graphs and variable selection with the {L}asso.
\newblock {\em Annals of Statistics}, 34:1436--1462, 2006.

\bibitem{NegahbanRWY10}
S.~Negahban, P.~Ravikumar, M.~J. Wainwright, and B.~Yu.
\newblock A unified framework for high-dimensional analysis of {M}-estimators
  with decomposable regularizers.
\newblock {\em Statistical Science}, 27:538--557, 2012.

\bibitem{RaskuttiWY11}
G.~Raskutti, M.~J. Wainwright, and B.~Yu.
\newblock Minimax rates of estimation for high--dimensional linear regression
  over $\ell_q$--balls.
\newblock {\em IEEE Trans. Info. Theory}, 57:6976--6994, 2011.

\bibitem{RavikumarWRY08}
P.~Ravikumar, M.~J. Wainwright, G.~Raskutti, and B.~Yu.
\newblock Model selection in gaussian graphical models: High-dimensional
  consistency of $\ell_1$-regularized {MLE}.
\newblock {\em In Advances in Neural Information Processing Systems (NIPS)},
  21, 2008.

\bibitem{RochaZY08}
G.~Rocha, P.~Zhao, and B.~Yu.
\newblock A path following algorithm for sparse pseudo-likelihood inverse
  covariance estimation (splice).
\newblock Technical report, University of California, Berkeley, 2008.

\bibitem{RothmanBLZ08}
A.~Rothman, P.~Bickel, E.~Levina, and J.~Zhu.
\newblock Sparse permutation invariant covariance estimation.
\newblock {\em Electronic Journal of Statistics}, 2:494--515, 2008.

\bibitem{SunZ11}
T.~Sun and C.-H. Zhang.
\newblock Scaled sparse linear regression.
\newblock {\em Biometrika}, 99:879--898, 2012.

\bibitem{vandeGeer07}
S.~van~de Geer.
\newblock The deterministic {L}asso.
\newblock Technical Report 140, ETH Zurich, Switzerland, 2007.

\bibitem{vandeGeerB09}
S.~van~de Geer and P.~B{\"u}hlmann.
\newblock On the conditions used to prove oracle results for the {L}asso.
\newblock {\em Electronic Journal of Statistics}, 3:1360--1392, 2009.

\bibitem{YangK10}
S.~Yang and E.~D. Kolaczyk.
\newblock Target detection via network filtering.
\newblock {\em IEEE Transactions on Information Theory}, 56(5):2502--2515,
  2010.

\bibitem{YeZ10}
F.~Ye and C.-H. Zhang.
\newblock Rate minimaxity of the{L}asso and {D}antzig selector for the $\ell_q$
  loss in $\ell_r$ balls.
\newblock {\em Journal of Machine Learning Research}, 11:3481--3502, 2010.

\bibitem{Yuan10}
M.~Yuan.
\newblock Sparse inverse covariance matrix estimation via linear programming.
\newblock {\em Journal of Machine Learning Research}, 11:2261--2286, 2010.

\bibitem{YuanL07}
M.~Yuan and Y.~Lin.
\newblock Model selection and estimation in the gaussian graphical model.
\newblock {\em Biometrika}, 94(1):19--35, 2007.

\bibitem{Zhang10}
C.-H. Zhang.
\newblock Nearly unbiased variable selection under minimax concave penalty.
\newblock {\em Annals of Statistics}, 38:894--942, 2010.

\bibitem{ZhangH08}
C.-H. Zhang and J.~Huang.
\newblock The sparsity and bias of the {L}asso selection in high-dimensional
  linear regression.
\newblock {\em Annals of Statistics}, 36(4):1567--1594, 2008.

\bibitem{ZhangZ11}
C.-H. Zhang and T.~Zhang.
\newblock A general theory of concave regularization for high dimensional
  sparse estimation problems.
\newblock {\em Statistical Science}, 27(4):576--593, 2012.

\bibitem{Zhang09-l1}
T.~Zhang.
\newblock Some sharp performance bounds for least squares regression with
  {$L_1$} regularization.
\newblock {\em Ann. Statist.}, 37(5A):2109--2144, 2009.

\end{thebibliography}
\end{document}